\documentclass[a4paper, 12pt]{article}
\usepackage{amsfonts}
\usepackage{amssymb}
\usepackage{amsthm}
\usepackage[dvips]{graphics,color}

\newtheorem{theorem}{Theorem}[section]
\newtheorem{lemma}[theorem]{Lemma}
\newtheorem{defn}[theorem]{Definition}
\newtheorem{cor}[theorem]{Corollary}
\newtheorem{prop}[theorem]{Proposition}
\theoremstyle{definition}{\newtheorem{eg}[theorem]{Example}}

\newcommand{\hf}[1][1]{\frac{#1}{2}}

\newcommand{\di}{\,\mathrm{d}}
\newcommand{\rmd}{\,\mathrm{d}}
\newcommand{\Or}{\mathrm{O}}
\newcommand{\rmi}{\mathrm{i}}

\begin{document}
\setlength{\unitlength}{10mm}
\title{ A one dimensional analysis of singularities and turbulence for the stochastic Burgers equation in $d$-dimensions}
\author{A D Neate $\qquad$  A Truman\\
\small{Department of Mathematics, University of Wales Swansea,}\\ [-0.8ex]
\small{Singleton Park, Swansea, SA2 8PP, Wales, UK.}}
\maketitle

\begin{abstract}

The inviscid limit of the stochastic Burgers equation, with body
forces white noise in time, is discussed in terms of the level
surfaces of the minimising Hamilton-Jacobi function, the classical
mechanical caustic and the Maxwell set and their algebraic
pre-images under the classical mechanical flow map. The problem is
analysed in terms of a reduced (one dimensional) action function. We
give an explicit expression for an algebraic surface containing the
Maxwell set and caustic in the polynomial case. Those parts of the
caustic and Maxwell set which are singular are characterised. We
demonstrate how the geometry of the caustic, level surfaces and
Maxwell set can change infinitely rapidly causing turbulent
behaviour which is stochastic in nature, and we determine its
intermittence in terms of the recurrent behaviour of two processes.

\end{abstract}

\section{Introduction}

Burgers equation has been used in studying turbulence and in
modelling the large scale structure of the universe \cite{Arnold4,
E, Shandarin}, as well as to obtain  detailed asymptotics for
stochastic Schr\"{o}dinger and heat equations
\cite{Elworthy,Elworthy2,Truman,Truman2,Truman3,Truman4}. It has
also played a part in Arnol'd's work on caustics and Maslov's works
in semiclassical quantum mechanics
\cite{Arnold,Arnold2,Maslov,Maslov2}.

   We consider the stochastic viscous Burgers equation for the
velocity field $v^{\mu}(x,t)\in\mathbb{R}^d$, where
$x\in\mathbb{R}^d$ and $t>0$,
\[\frac{\partial v^{\mu}}{\partial
t}+\left(v^{\mu}\cdot \nabla\right) v^{\mu} = \hf[\mu^2]\Delta
v^{\mu} - \nabla V(x) -\epsilon \nabla k_t(x)\dot{W}_t, \qquad
v^{\mu}(x,0) = \nabla S_0 (x)+\Or(\mu^2).\] Here $\dot{W}_t$ denotes
white noise and $\mu^2$ is the coefficient of viscosity which we
assume to  be small. We are interested in the advent of
discontinuities in the inviscid limit of the Burgers fluid velocity
$v^0(x,t)$  where  $v^{\mu}(x,t)\rightarrow v^0(x,t)$ as
$\mu\rightarrow 0.$

Using the Hopf-Cole transformation $v^{\mu}(x,t)=-\mu^2\nabla\ln
u^{\mu}(x,t)$, the Burgers equation becomes the Stratonovich heat
equation,
\[\frac{\partial u^{\mu}}{\partial
t}=\hf[\mu^2]\Delta u^{\mu}
+\mu^{-2}V(x)u^{\mu}+\frac{\epsilon}{\mu^{2}}
k_t(x)u^{\mu}\circ\dot{W}_t , \quad u^{\mu}(x,0) =
\exp\left(-\frac{S_0(x)}{\mu^2}\right)T_0(x),\] where the
convergence factor $T_0$ is related to the initial Burgers fluid
density \cite{Hopf}.

Now let, \begin{equation}\label{action} A[X]: =
\hf\int_0^t \dot{X}^2(s)\di s -\int_0^t V(X(s)) \di s -\epsilon
\int_0^t k_s(X(s))\di W_s,\end{equation}
and select a path $X$ which minimises $A[X]$. This requires,
\begin{equation}\label{ie3}
\di \dot{X}(s)+\nabla V(X(s)) \di s
+\epsilon \nabla k_s(X(s))\di W_s=0.
\end{equation}
We then define the stochastic action
$A(X(0),x,t):=\inf\limits_X\left\{ A[X]:X(t)=x\right\}.$
Setting,
\[\mathcal{A}(X(0),x,t):=S_0(X(0))+A(X(0),x,t),\]
and then minimising $\mathcal{A}$  over $X(0)$, gives
$\dot{X}(0) = \nabla
S_0(X(0)).$
Moreover, it follows that, \[\mathcal{S}_t(x):=\inf\limits_{X(0)}
\left\{\mathcal{A}(X(0),x,t)\right\},\]
is  the minimal solution of the Hamilton-Jacobi
equation,
\begin{equation}\label{ie4}
\di \mathcal{S}_t +\left(\hf|\nabla \mathcal{S}_t|^2+V(x)\right)\di t
+\epsilon k_t(x) \di W_t=0,
\qquad \mathcal{S}_{t=0}(x) = S_0(x).\end{equation}
Following the work of Donsker, Freidlin et al \cite{Freidlin},
$-\mu^2\ln u^{\mu}(x,t) \rightarrow
\mathcal{S}_t(x)$ as $\mu\rightarrow 0$. This gives
the inviscid limit of the minimal entropy solution of Burgers equation  as
$v^0(x,t)=\nabla\mathcal{S}_t(x)$ \cite{Dafermos}.

\begin{defn}
The stochastic wavefront at time $t$ is defined to be the set,
\[\mathcal{W}_t=\left\{x:\quad\mathcal{S}_t(x)=0\right\}.\]
\end{defn}
For small $\mu$ and fixed $t$,  $u^{\mu}(x,t)$ switches
continuously from being exponentially large to small as $x$ crosses the
wavefront $\mathcal{W}_{t}$. However, $u^{\mu}$ and $v^{\mu}$ can
also switch discontinuously.

Define the classical flow map
$\Phi_s:\mathbb{R}^d\rightarrow\mathbb{R}^d$ by,
\[\rmd \dot{\Phi}_s +\nabla V(\Phi_s) \rmd s+\epsilon\nabla k_s(\Phi_s)\rmd
W_s =0,\qquad\Phi_0 = \mbox{id},\qquad \dot{\Phi}_0 = \nabla S_0.\]
Since $X(t) =x$ it follows that $X(s) = \Phi_s\left(
\Phi_t^{-1}( x)\right)\!,$ where the pre-image $x_0(x,t) = \Phi_t^{-1} (x)$ is
not necessarily unique.

Given some regularity and boundedness, the global inverse function
theorem gives a caustic time $T(\omega)$ such that for
$0<t<T(\omega)$, $\Phi_t$ is a random diffeomorphism; before the
caustic time $v^0(x,t) = \dot{\Phi}_t\left(\Phi_t^{-1}(x)\right)$ is
the inviscid limit of a classical solution of the Burgers equation
with probability one.

The method of characteristics suggests that discontinuities in
$v^0(x,t)$ are associated with the non-uniqueness of the real pre-image $x_0(x,t)$.
When this occurs,  the classical flow map
$\Phi_t$ focusses
an infinitesimal volume of points $\rmd
x_0$ into a zero volume $\rmd X(t)$.
 \begin{defn}\label{i2} The caustic at time $t$ is defined to be the set,
\[ C_t = \left\{ x: \quad\det\left(\frac{\partial X(t)}{\partial x_0}\right) = 0\right\}. \]
\end{defn}

 Assume that $x$ has $n$ real pre-images,
\[\Phi_t^{-1}\left\{x\right\} =
\left\{x_0(1)(x,t),x_0(2)(x,t),\ldots,x_0(n)(x,t)\right\},\]  where
each $x_0(i)(x,t)\in\mathbb{R}^d$. Then the Feynman-Kac formula and
Laplace's method in infinite dimensions give for a non-degenerate
critical point, \begin{equation}\label{useries}u^{\mu}(x,t)=
\sum\limits_{i=1}^n \theta_i
\exp\left(-\frac{S_0^i(x,t)}{\mu^2}\right),
\end{equation} where
$S_0^i(x,t) :=
S_0\left(x_0(i)(x,t)\right)+A\left(x_0(i)(x,t),x,t\right),$
and $\theta_i$ is an asymptotic series in $\mu^2$. An asymptotic
series in $\mu^2$  can also be found for
$v^{\mu}(x,t)$ \cite{Truman5}.
Note that $\mathcal{S}_t(x) = \min
\{S_0^i(x,t):i=1,2,\ldots,n\}$.

\begin{defn}\label{i3}
The Hamilton-Jacobi level surface is the set,
\[H_t^c = \left\{ x:\quad S_0^i(x,t) =c \mbox{ for
some }i\right\}.\]
The zero level surface $H_t^0$ includes
the wavefront $\mathcal{W}_t$.
\end{defn}
As $\mu\rightarrow 0$, the dominant term in the expansion (\ref{useries}) comes
from the minimising $x_0(i)(x,t)$ which we denote $\tilde{x}_0(x,t)$.  Assuming $\tilde{x}_0(x,t)$ is unique, we obtain the inviscid limit of
the Burgers  fluid velocity as
$v^0(x,t) = \dot{\Phi}_t\left(\tilde{x}_0(x,t)\right).$

If the minimising pre-image $\tilde{x}_0(x,t)$ suddenly changes value between two pre-images $x_0(i)(x,t)$ and $x_0(j)(x,t)$, a jump discontinuity will also occur in the inviscid limit of the Burgers fluid velocity.
There are  two distinct ways in which the minimiser can change; either two pre-images coalesce and disappear (become complex), or the minimiser switches between two pre-images at the same action value. The first of these occurs as $x$ crosses the caustic and when the minimiser disappears the caustic is said to be cool. The second occurs as $x$ crosses the Maxwell set and again, when the minimiser is involved the Maxwell set is said to be cool.

 \begin{defn}\label{m1}
The Maxwell set is given by,
\begin{eqnarray*}
M_t & = & \left\{x:\, \exists \,x_0,\check{x}_0\in\mathbb{R}^d \mbox{ s.t. }\right.\\
&   &\quad\left.
x=\Phi_t(x_0)=\Phi_t(\check{x}_0), \,x_0\neq \check{x}_0 \mbox{ and }
 \mathcal{A}(x_0,x,t)=\mathcal{A}(\check{x}_0,x,t)
\right\}.
\end{eqnarray*}
 \end{defn}

\begin{eg}[The generic Cusp] Let $V(x,y)=0,$ $k_t(x,y)=0$ and $S_0(x_0,y_0)=x_0^2 y_0/2$.
This initial condition leads to the \emph{generic Cusp}, a semicubical parabolic caustic shown in Figure \ref{cusp}. The caustic $C_t$ (long dash) is given by,
\[
 x_t(x_0)  =  t^2 x_0^3,\quad
 y_t(x_0)  = \hf[3]tx_0^2-\frac{1}{t}.
\]
 The zero level surface $H_t^0$ (solid line) is,
\[
 x_{(t,0)}(x_0)   =  \hf[x_0]\left(1\pm\sqrt{1-t^2 x_0^2}\right),\quad
 y_{(t,0)}(x_0)  =  \frac{1}{2t}\left(t^2x_0^2-1\pm
\sqrt{1-t^2x_0^2}\right),\]
and the Maxwell set $M_t$ (short dash) is $x=0$ for  $y>-1/t.$

\begin{figure}[h!]
\begin{center}\begin{picture}(5.5,5.5)
    \put(0,0){\resizebox{5.5cm}{!}{\includegraphics{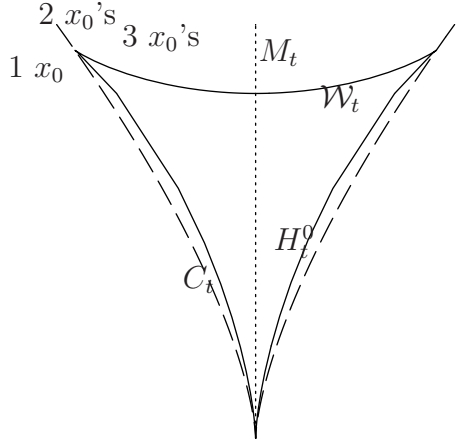}}}
    \put(1.8,2){$C_t$}
    \put(-0.5,4.8){$1$ $x_0$}
    \put(1,5.2){$3$ $x_0$'s}
    \put(-0.1,5.5){$2$ $x_0$'s}
    \put(3.6,4.4){$\mathcal{W}_t$}
    \put(3,2.5){$H_t^0$}
    \put(2.8,5){$M_t$}
  \end{picture}
  \end{center}
\caption{Cusp and Tricorn.}\label{cusp}
\end{figure}
\end{eg}

\noindent\emph{Notation:} Throughout this paper $x,x_0,x_t$ etc will
denote vectors, where normally $x=\Phi_t (x_0)$. Cartesian
coordinates of these will be indicated using a sub/superscript where
relevant; thus $x=(x_1,x_2,\ldots,x_d)$,
$x_0=(x_0^1,x_0^2,\ldots,x_0^d)$ etc. The only exception  will be
 in discussions of explicit examples in two and three dimensions when
we will use $(x,y)$ and $(x_0,y_0)$ etc to denote the vectors.

\section{Some background}

We begin by summarising  some of the geometrical results established
by Davies, Truman and Zhao (DTZ) \cite{Davies,Davies2,Davies3} and
presenting some minor generalisations of their results
\cite{Neate2,Neate3}. Following equation (\ref{action}), let the
stochastic action be defined,
\[
A(x_0,p_0,t) =  \hf\int_0^t \dot{X}(s)^2\rmd s
 -\int_0^t\Bigg[ V(X(s))\rmd s +\epsilon k_s(X(s))\rmd
W_s\Bigg],\]
where
$X(s)=X(s,x_0,p_0)\in\mathbb{R}^d$ and,
\[\rmd\dot{X}(s) = -\nabla V(X(s))\rmd s
-\epsilon \nabla k_s(X(s))\rmd W_s,\quad
X(0)=x_0, \quad\dot{X}(0)=p_0,\]
for $ s\in[0,t]$ with  $x_0,p_0\in\mathbb{R}^d$. We assume
$X(s)$ is $\mathcal{F}_s$ measurable and unique.

\begin{lemma} \label{i8}Assume $S_0,V\in
C^2$ and $k_t\in C^{2,0}$, $\nabla V,\nabla k_t$ Lipschitz with Hessians
$\nabla^2 V,\nabla^2 k_t$ and all second derivatives with respect to
space variables of $V$ and $k_t$ bounded. Then for $p_0$, possibly
$x_0$ dependent,
\[\frac{\partial A}{\partial
x_0^{\alpha}}(x_0,p_0,t) =\dot{X}(t)\cdot \frac{\partial
X(t)}{\partial x_0^{\alpha}}-\dot{X}_{\alpha}(0),\qquad\alpha
=1,2,\ldots,d.\]
\end{lemma}

 Methods of Kolokoltsov et al
\cite{Kolokoltsov,Kunita}  guarantee that for small $t$ the map $p_0\mapsto
X(t,x_0,p_0)$ is onto for all $x_0$. Therefore, we can define,
\[A(x_0,x,t)
=\left.A(x_0,p_0,t)\right|_{p_0=p_0(x_0,x,t){\displaystyle,}}\] where
$p_0=p_0(x_0,x,t)$ is the random minimiser (assumed unique) of
$A(x_0,p_0,t)$ when $X(t,x_0,p_0)=x$. The
stochastic action corresponding to the initial momentum $\nabla
S_0(x_0)$ is then $\mathcal{A}(x_0,x,t) := A(x_0,x,t)+S_0(x_0).$

\begin{theorem}\label{sflow} If $\Phi_t$ is the stochastic flow map, then $\Phi_t(x_0)=x$ is equivalent to,
\[\frac{\partial}{\partial x_0^{\alpha}}
\left[\mathcal{A}(x_0,x,t)\right]=0,\qquad \alpha=1,2,\ldots,d.\]
\end{theorem}

The Hamilton-Jacobi level surface $H_t^c$ is
obtained by eliminating $x_0$ between,
\[\mathcal{A}(x_0,x,t) =c \quad \mbox{and}\quad
\frac{\partial\mathcal{A}}{\partial x_0^{\alpha}}(x_0,x,t) =0, \quad
\alpha=1,2,\ldots, d.\] Alternatively, if we eliminate $x$ to give
an expression in $x_0$, we have the pre-level surface
$\Phi_t^{-1}H_t^c$. Similarly the caustic $C_t$ (and pre-caustic
$\Phi_t^{-1}C_t$) are obtained by eliminating $x_0$ (or $x$) between,
\[\det\left(
\frac{\partial^2\mathcal{A}} {\partial x_0^{\alpha}\partial
x_0^{\beta}}(x_0,x,t) \right)_{\alpha,\beta=1,2,\ldots, d} =0 \quad
\mbox{and} \quad \frac{\partial\mathcal{A}}{\partial
x_0^{\alpha}}(x_0,x,t) =0\quad \alpha = 1,2,\ldots, d.\] These
pre-images are calculated algebraically which are not necessarily the
topological inverse images of the surfaces $C_t$ and $H_t^c$ under
$\Phi_t$.

Assume that $A(x_0,x,t)$ is $C^4$ in space variables with
$\det\left( \frac{\partial^2\mathcal{A}}{\partial
x_0^{\alpha}\partial x^{\beta}}\right)\neq 0.$

\begin{defn}
A curve $x=x(\gamma)$, $\gamma \in N(\gamma_0,\delta)$, is said to
have a generalised cusp at $\gamma= \gamma_0$, $\gamma$ being an
intrinsic variable such as arc length, if
$\frac{\di x}{\di\gamma}( \gamma_0)=0.$
\end{defn}
\begin{lemma}\label{s11}
Let $\Phi_t$ denote the flow map and let  $\Phi_t^{-1}\Gamma_t$ and
$\Gamma_t$ be some surfaces where if $x_0\in\Phi_t^{-1}\Gamma_t$
then $x=\Phi_t (x_0)\in\Gamma_t$. Then $\Phi_t$ is a differentiable
map from $\Phi_t^{-1}\Gamma_t$ to $\Gamma_t$ with Frechet derivative,
\[D\Phi_t(x_0) = \left(-\frac{\partial^2\mathcal{A}}{\partial x\partial x_0}(x_0,x,t)\right)^{-1}
\left(\frac{\partial^2\mathcal{A}}{\partial x_0^2}(x_0,x,t)\right).
\]
\end{lemma}
\begin{lemma}
Let $x_0(s)$ be any two dimensional intrinsically parameterised
curve, and define $x(s)=\Phi_t(x_0(s)).$
Let $e_0$ denote the zero eigenvector of
$\left(\frac{\partial^2\mathcal{A}}{(\partial x_0)^2}\right) $ and
assume that $\ker \left(\frac{\partial^2\mathcal{A}}{(\partial
x_0)^2}\right) =\langle e_0 \rangle$. Then, there is a generalised
cusp on $x(s)$ when $s=\sigma$ if and only if either:
\begin{enumerate}
\item there is a generalised cusp on $x_0(s)$ when $s=\sigma$; or,
\item $x_0(\sigma)$ is on the pre-caustic and the tangent
$\frac{\di x_0}{\di s}(s)$ at $s=\sigma$ is parallel to $e_0$.
\end{enumerate}
\end{lemma}

\begin{prop}\label{i10}The normal to $\Phi_t^{-1}H_t^c$ is,
\[n(x_0) =
-\left(\frac{\partial^2\mathcal{A}}{\partial x_0\partial x_0}\right)
\left(\frac{\partial^2\mathcal{A}}{\partial x_0\partial
x}\right)^{-1} \dot{X}\left(t,x_0,\nabla S_0(x_0)\right).\]
\end{prop}
\begin{cor}\label{i11a}
In two dimensions, let $\Phi_t^{-1}H_t^c$ meet $\Phi_t^{-1}C_t$ at  $x_0$ where $n(x_0)\neq 0$ and $\ker\left(\frac{\partial^2\mathcal{A}}{(\partial x_0)^2}\right)=\langle e_0\rangle$. Then the tangent to $\Phi_t^{-1}H_t^c$ at $x_0$  is parallel to $e_0$.\end{cor}

\begin{prop}\label{i11}
In two dimensions, assume that $n(x_0)\neq 0$ where $x_0\in\Phi_t^{-1} H_t^c$, so that $\Phi_t^{-1}H_t^c$ does not have a generalised cusp at $x_0$. Then $H_t^c$ can only have a generalised cusp at $\Phi_t(x_0)$ if $\Phi_t(x_0)\in C_t$. Moreover, if $x=\Phi_t(x_0)\in\Phi_t\left\{\Phi_t^{-1}C_t\cap H_t^{-1}\right\}$ then $H^c_t$ will have a generalised cusp.
\end{prop}

\begin{figure}[h!]
 \begin{center}
\begin{tabular}{cc}
 \begin{picture}(5.5,5.5)
 \put(0,0){\resizebox{55mm}{!}{\includegraphics{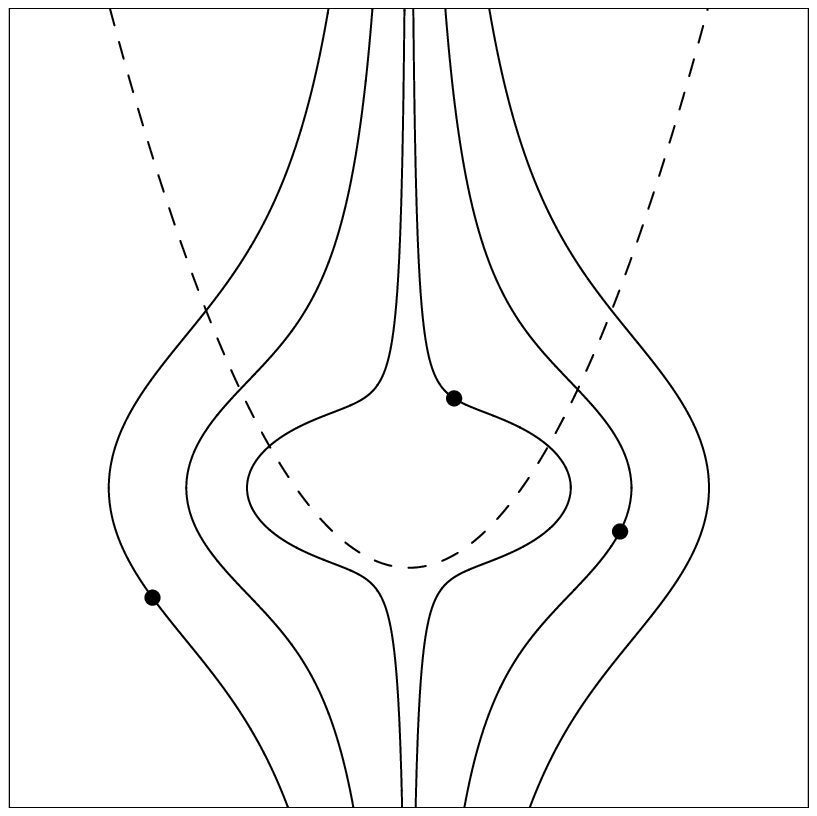}}}
 \put(0.1,0.3){(a)}
 \end{picture}&
\begin{picture}(5.5,5.5)
 \put(0,0){\resizebox{55mm}{!}{\includegraphics{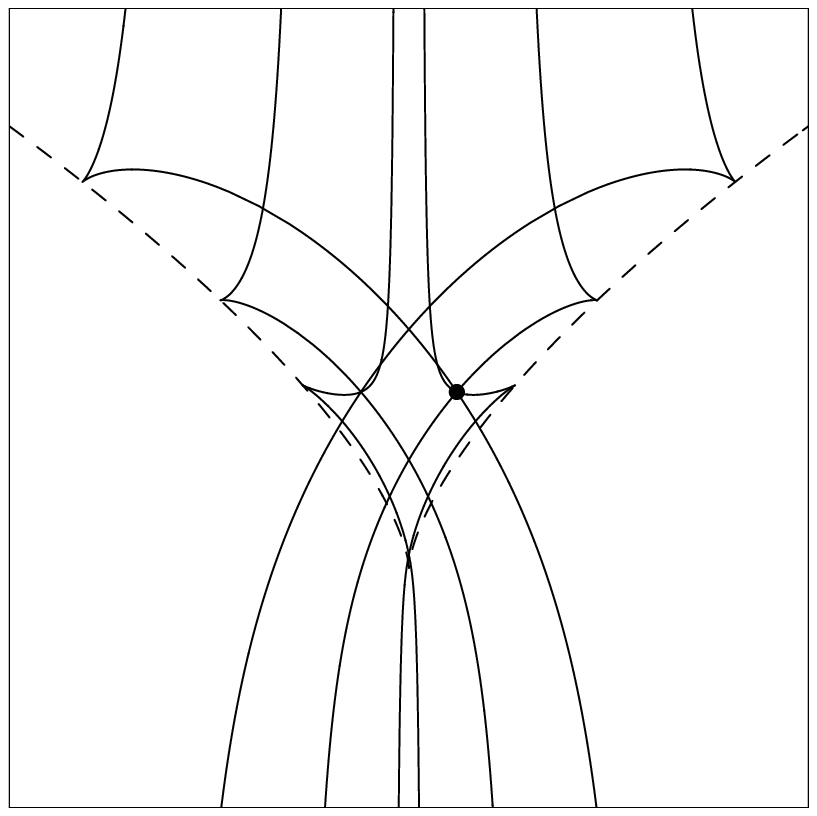}}}
 \put(0.1,0.3){(b)}
 \end{picture}
\end{tabular}
\end{center}
 \caption{(a) The pre-level surface (solid line) and pre-caustic
 (dashed),
 (b) the level surface (solid line) and caustic (dashed), both for the generic Cusp with $c>0$.}
 \label{if2}
 \end{figure}

\begin{eg}[The generic Cusp]
Figure \ref{if2} shows that a point lying on three level surfaces  has three distinct real pre-images each on a separate pre-level surface. A cusp only occurs on the corresponding level surface when the pre-level surface intersects the pre-caustic. Thus, a level surface only has a
cusp on the caustic, but it does not have to be cusped when it meets the caustic.
\end{eg}

\begin{theorem}Let,
\[
x \in  \mbox{Cusp}\left(H_t^c\right)
=
\left\{x\in\Phi_t\left(\Phi_t^{-1}C_t\cap\Phi_t^{-1}H^c_t\right),x=\Phi_t(x_0),
n(x_0)\neq 0\right\}.\]Then in three dimensions in the
stochastic case, with probability one, $T_x$ the tangent space to
the level surface at $x$ is at most one dimensional.
\end{theorem}

\section{A one dimensional analysis}
In this section we  outline a one dimensional analysis first described by
Reynolds, Truman and Williams (RTW) \cite{Truman6}. \begin{defn}\label{i16}
    The $d$-dimensional flow map $\Phi_{t}$ is globally reducible if for
    any
    $x=(x_{1},x_{2},\ldots ,x_{d})$ and $x_{0}=(x_{0}^{1},x_{0}^{2},\ldots,
    x_{0}^{d})$ where $x=\Phi_{t}(x_{0}),$ it is possible to write
    each coordinate $x_0^\alpha$ as a function of the lower
    coordinates. That is,
    \begin{equation}\label{ie16}
    x=\Phi_{t}(x_{0}) \quad \Rightarrow\quad
        x_{0}^{\alpha}=x_{0}^{\alpha}(x,x_{0}^{1},x_{0}^{2},\ldots x_{0}^{\alpha-1},t)\mbox{ for $ \alpha=d,d-1,\ldots ,2.$}
    \end{equation}

\end{defn}
Therefore, using Theorem \ref{sflow}, the flow map is globally
reducible if we can find a chain of $C^{2}$ functions
$x_{0}^{d},x_{0}^{d-1},\ldots,x_0^2$ such that,
\begin{eqnarray*}
     x_{0}^{d}=x_{0}^{d}(x,x_{0}^{1},x_{0}^{2},\ldots
   x_{0}^{d-1},t)\quad
&\Leftrightarrow &\quad
\frac{\partial\mathcal{A}}{\partial x_{0}^{d}}(x_{0},x,t)=0,\\
 x_{0}^{d-1}=x_{0}^{d-1}(x,x_{0}^{1},x_{0}^{2},\ldots
x_{0}^{d-2},t) \quad&\Leftrightarrow &\quad
\frac{\partial\mathcal{A}}{\partial
x_{0}^{d-1}}(x_{0}^{1},x_{0}^{2},\ldots,x_{0}^{d}(\ldots),x,t)=0,\\
& \vdots &\\
 x_{0}^{2}=x_{0}^{2}(x,x_{0}^{1},t)\quad &\Leftrightarrow &
 \end{eqnarray*}
\[\qquad\qquad\qquad\qquad
\qquad \qquad \qquad \,\,\,\,\,\frac{\partial\mathcal{A}}{\partial
x_{0}^{2}}(x_{0}^{1},x_{0}^{2},x_{0}^{3}(x,x_{0}^{1},x_{0}^{2},t),\ldots,x_{0}^{d}(\ldots),x,t)=0,
\]
where $x_{0}^{d}(\ldots)$ is the expression only involving $x_0^1$
and $x_0^2$ gained by substituting each of the functions
$x_0^3,\ldots,x_0^{d-1}$ repeatedly into
$x_{0}^{d}(x,x_{0}^{1},x_{0}^{2},\ldots, x_{0}^{d-1},t)$. This
requires that no roots are repeated to ensure that none of
the second derivatives of $\mathcal{A}$ vanish. We assume also
that there is a favoured ordering of coordinates and a corresponding
decomposition of $\Phi_{t}$ which allows the non-uniqueness to be
reduced to the level of the $x_{0}^{1}$ coordinate.
This assumption   appears to be quite restrictive. However,
local reducibility at $x$ follows from the implicit function theorem and some mild assumptions on the derivatives of $\mathcal{A}$.

\begin{defn}\label{i17}
If $\Phi_{t}$ is globally reducible then  the reduced
    action function is the univariate function gained from evaluating the
     action with equations (\ref{ie16}),
    \[f_{(x,t)}(x_0^1):=f(x_{0}^{1},x,t) =\mathcal{A}
    (x_{0}^{1},x_{0}^{2}(x,x_{0}^{1},t),x_{0}^{3}(\ldots),\ldots,x,t).\]
\end{defn}

\begin{lemma} \label{i18}If $\Phi_t$ is globally reducible, modulo the above assumptions,
\begin{eqnarray*}\lefteqn{\left|
    \det
    \left.\left(\frac{\partial^{2}\mathcal{A}}{(\partial
    x_{0})^{2}}({x}_{0},{x},t)\right)
    \right|_{x_0=(x_0^1,x_0^2(x,x_0^1,t),\ldots,x_0^d(\ldots))}\right|}\\
    &  = & \prod\limits_{\alpha=1}^{d}\left|\left[
    \left(\frac{\partial}{\partial
 x_{0}^{\alpha}}\right)^2\!\!\! \mathcal{A}
 (x_{0}^{1},\ldots,x_{0}^{\alpha},x_0^{\alpha+1}(\ldots),\ldots
    ,x_{0}^{d}(\ldots),{x},t)\right]_{\begin{array}{c}
   { \scriptstyle x_0^2=x_0^2(x,x_0^1,t)}\\[-1ex]
     {\scriptstyle{\vdots}}\\[-1ex]
   {\scriptstyle x_0^{\alpha} = x_0^{\alpha}(\ldots)}
    \end{array}}\right|
\end{eqnarray*}
   where the first term is $f_{({x},t)}''(x_{0}^{1})$ and the last
    $d-1$ terms are non zero.
\end{lemma}

\begin{theorem}\label{i19}
Let the classical mechanical flow map $\Phi_t$ be globally
reducible.   Then:
    \begin{enumerate}
    \item $ f_{(x,t)}' (x_{0}^{1})=0$
     and the  equations (\ref{ie16})
     \emph{$\Leftrightarrow x=\Phi_{t} (x_{0}),$}
    \item
    $ f_{(x,t)}' (x_{0}^{1})
    = f_{(x,t)}'' (x_{0}^{1})=0$
    and the  equations (\ref{ie16})\\
    \emph{$\Leftrightarrow x=\Phi_{t}(x_{0})$} is such that the number of real
    solutions $x_{0}$ changes.
    \end{enumerate}
\end{theorem}

\section{Analysis of the caustic}
We begin by parameterising the caustic
$0  =  \det\left(
D\Phi_{t}(x_0)
\right)
$
  from Definintion \ref{i2}; this equation only involves $x_0$ and $t$, and
is therefore the pre-caustic. We use this to parameterise the pre-caustic
 as,
\[x_0^1=\lambda_1,\quad x_0^2=\lambda_2,\quad\ldots,\quad x_0^{d-1}=\lambda_{d-1}\quad\mbox{and}\quad
x_0^d=x_{0}^d\left(\lambda_1,\lambda_2,\ldots,\lambda_{d-1}\right).\]
The parameters are restricted to be real so that only
 real pre-images are considered.

\begin{defn}\label{c1}
 For any
$\lambda=\left(\lambda_1,\lambda_2,\ldots,\lambda_{d-1}\right)\in\mathbb{R}^{d-1}$
the pre-parameterisation of the caustic is given by $
x_t(\lambda):=\Phi_t\left(\lambda,x_{0}^d(\lambda)\right).$
\end{defn}
The  pre-parameterisation will be intrinsic if  ker$(D\Phi_t)$ is one dimensional.

\begin{cor}\label{c2}
    Let $x_{t}(\lambda)$ denote the pre-parameterisation of the caustic where $\lambda
    =(\lambda_{1},\lambda_{2},\ldots,\lambda_{d-1})\in\mathbb{R}^{d-1}$. Then $f'_{(x_{t}(\lambda),t)}(\lambda_{1}) =
    f''_{(x_{t}(\lambda),t)}(\lambda_{1})=0.$
\end{cor}

\begin{prop}\label{c14}
Let $x_{t}(\lambda)$
    denote the pre-parameterisation of the caustic where $\lambda
    =(\lambda_{1},\lambda_{2},\ldots,\lambda_{d-1})\in\mathbb{R}^{d-1}$.
    Assume  $f_{(x_{t}(\lambda),t)}(x_{0}^{1})\in C^{p+1}$ then,
    in $d$-dimensions, if the tangent to the caustic is
    at most $(d-p+1)$-dimensional at $x_{t}(\tilde{\lambda})$,
    \[f'_{(x_{t}(\tilde{\lambda}),t)}(\tilde{\lambda}_{1}) =
    f''_{(x_{t}(\tilde{\lambda}),t)}(\tilde{\lambda}_{1}) =
    \ldots =
    f^{(p)}_{(x_{t}(\tilde{\lambda}),t)}(\tilde{\lambda}_{1}) =0.\]
\end{prop}
\begin{proof}
Follows by repeatedly differentiating
$f'''_{(x_t(\lambda),t)}(\lambda_1)=0,$ which holds if the tangent
space at $x_t(\lambda)$ is $(d-2)$-dimensional \cite{Neate2}.
\end{proof}

\begin{figure}[h!]
\setlength{\unitlength}{7.5mm}
\noindent\begin{center}\begin{tabular}{ccc} \textbf{Before Caustic}& \textbf{On
Cool Caustic} & \textbf{Beyond Caustic}\\
\begin{picture}(4.5,3) \put(0,0){\resizebox{33.75mm}{!}{\includegraphics{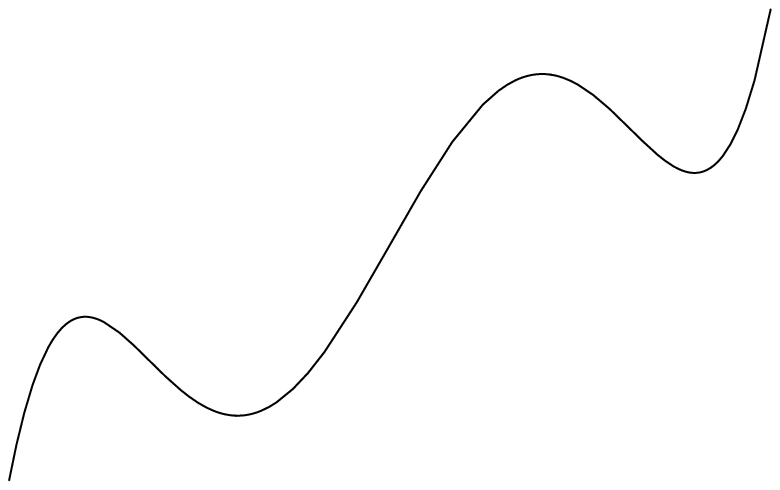}}}
\put(0.3,1.1){\scriptsize$x_0^1(1)$} \put(0.5,0.2){\scriptsize
$x_0^1(2)=\tilde{x}^1_0(x,t)$} \put(2.65,2.5){\scriptsize
$x_0^1(3)$} \put(3.7,1.5){\scriptsize $x_0^1(4)$}
\put(1.25,0.6){\scriptsize $(a)$} \put(3.8,2){\scriptsize $(b)$}
\end{picture}&
\begin{picture}(4.5,3) \put(0,0){\resizebox{33.75mm}{!}{\includegraphics{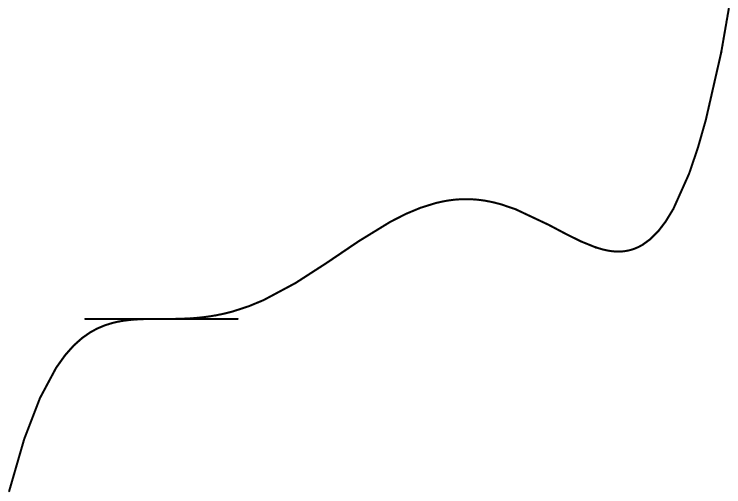}}}
\put(0.8,1.2){\scriptsize $(a)$}
\end{picture}&
\begin{picture}(4.5,3)
\put(0,0){\resizebox{33.75mm}{!}{\includegraphics{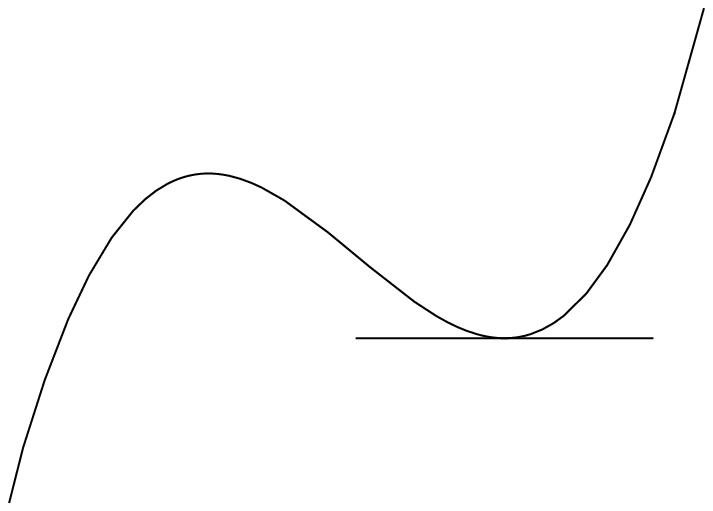}}}
\put(2.7,1.1){\scriptsize $(b)$}
\end{picture}\\
\emph{ Minimiser at }   & \emph{ Two $x_0^1$'s coalescing} &
\emph{ Minimiser jumps.} \\
      $x_0^1(2)(x,t)=\tilde{x}^1_0(x,t)$.              &\emph { form point of inflexion.} &
\end{tabular}\end{center}
 \caption{The graph of $f_{(x,t)}(x_0^1)$ as $x$ crosses the caustic.}\label{cf1}
\end{figure}

From Corollary \ref{c2}, there is a critical point of inflexion on
$f_{(x,t)}(x_0^1)$ at $x_{0}^{1}=\lambda_{1}$ when $x=x_t(\lambda)$.
Consider an example where for $x$ on one side of the caustic there
are  four real critical points on $f_{(x,t)}(x_0^1)=0$. Let them be
enumerated $x_0^1(i)(x,t)$ for $i=1$ to $4$ and denote the
minimising critical point $\tilde{x}_0^1(x,t)$. Figure \ref{cf1}
illustrates how the minimiser jumps from $(a)$ to $(b)$ as $x$
crosses the caustic.
This will cause $u^{\mu}$
and $v^{\mu}$ to jump for small $\mu$ and the caustic at such a
point is described as being cool.

\begin{defn}\label{c19}
Let   $x_{t}(\lambda)$ be the pre-parameterisation of the caustic. Then $x_t(\lambda)$ is on the cool part of the caustic if
    $f_{(x_t(\lambda),t)}(\lambda_1)\le f_{(x_t(\lambda),t)}(x_0^1(i)(x_t(\lambda),t))$
    for all $i=1,2,\ldots,n$ where $x_0^1(i)(x,t)$ denotes an enumeration of all the real roots for $x_0^1$ to
    $f_{(x,t)}'(x_0^1)=0.$ If  the caustic is not cool it is hot.
\end{defn}

\begin{figure}[h!]
\begin{tabular}[t]{c|cc|cc|}
\emph{Increasing $\lambda$} &\multicolumn{2}{c|}{\emph{Caustic
changes hot to cool}}&
\multicolumn{2}{c|}{\emph{No change in caustic}}\\
\begin{picture}(1,1.5)
\put(0.5,1.5){\vector(0,-1){1.5}}
\end{picture}
&\resizebox{22mm}{!}{\includegraphics{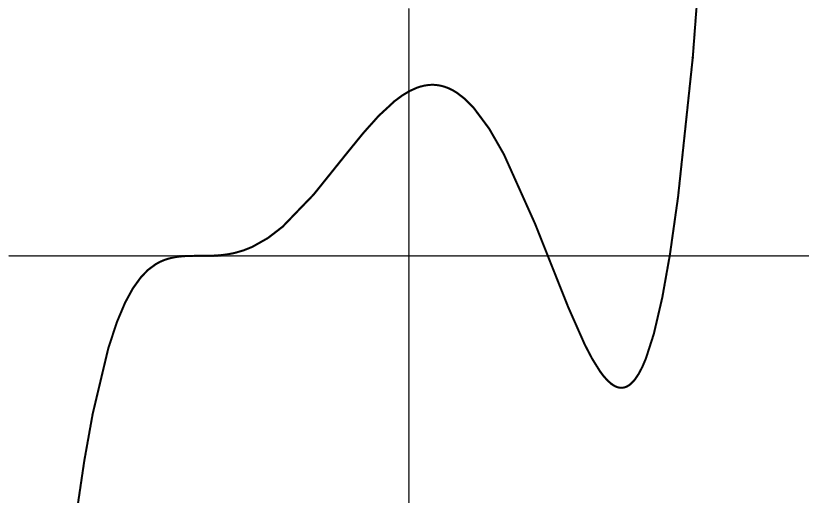}}&
 \resizebox{22mm}{!}{\includegraphics{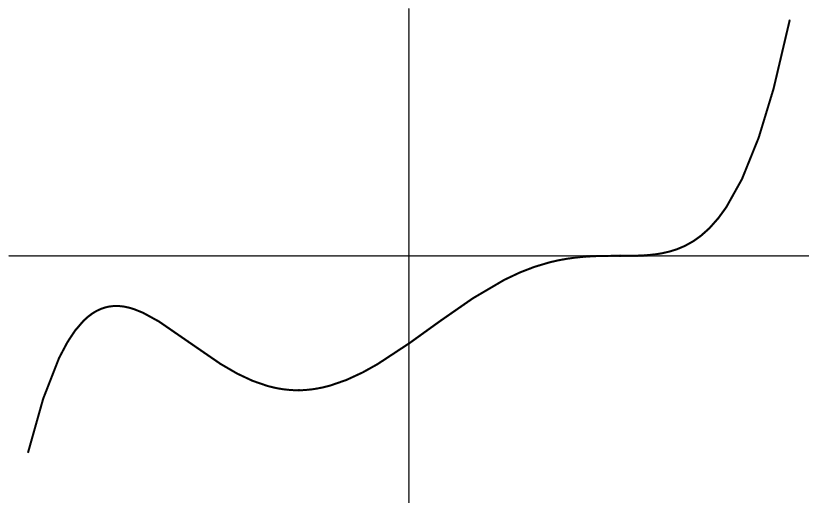}}&
 \resizebox{22mm}{!}{\includegraphics{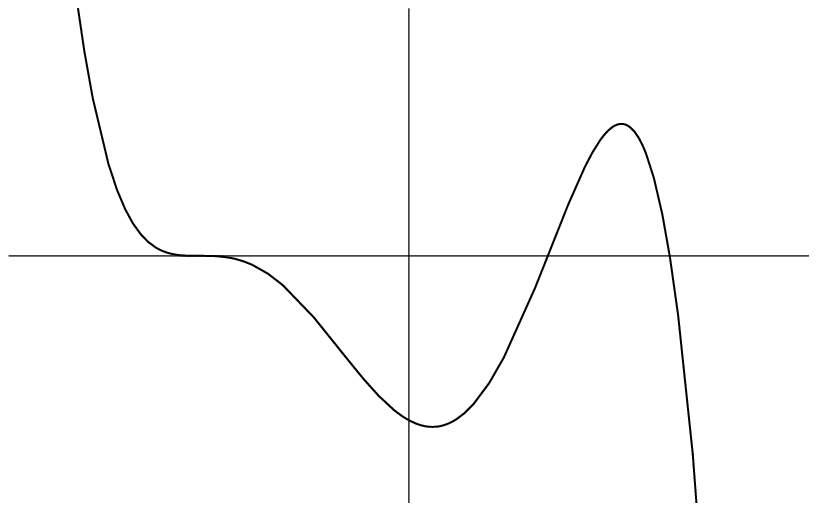}}&
 \resizebox{22mm}{!}{\includegraphics{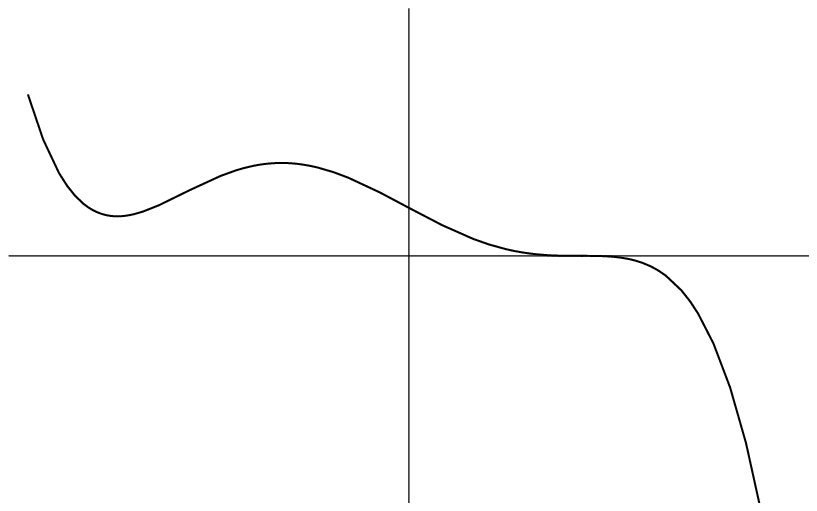}}
  \\
\raisebox{4ex}{\begin{tabular}{c}
\emph{Possible}\\
\emph{ hot/cool }\\
\emph{boundary}
\end{tabular}}
&  \resizebox{22mm}{!}{\includegraphics{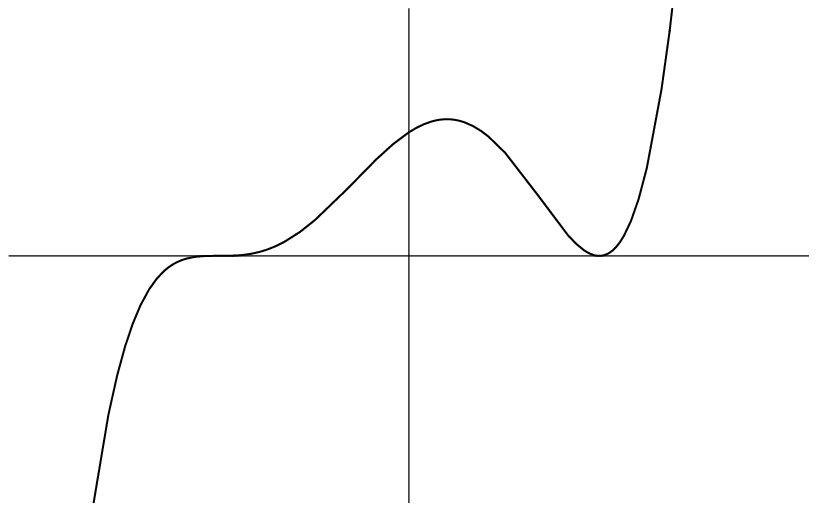}}&
 \resizebox{22mm}{!}{\includegraphics{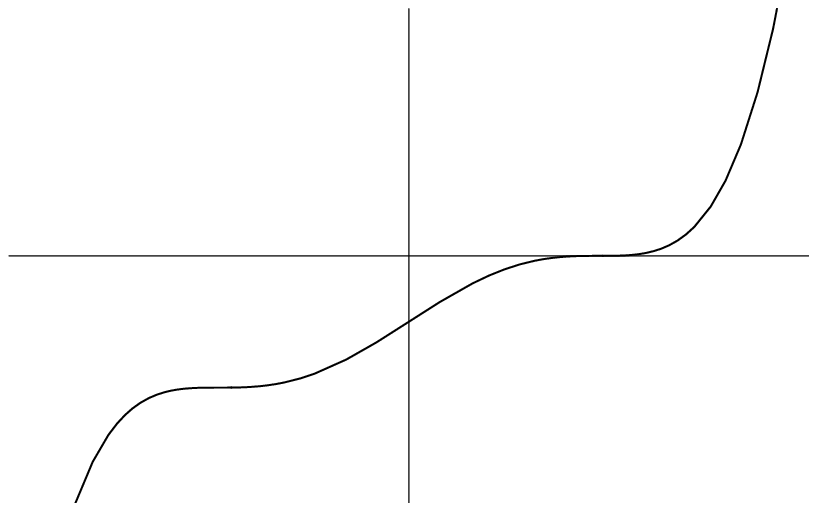}}&
 \resizebox{22mm}{!}{\includegraphics{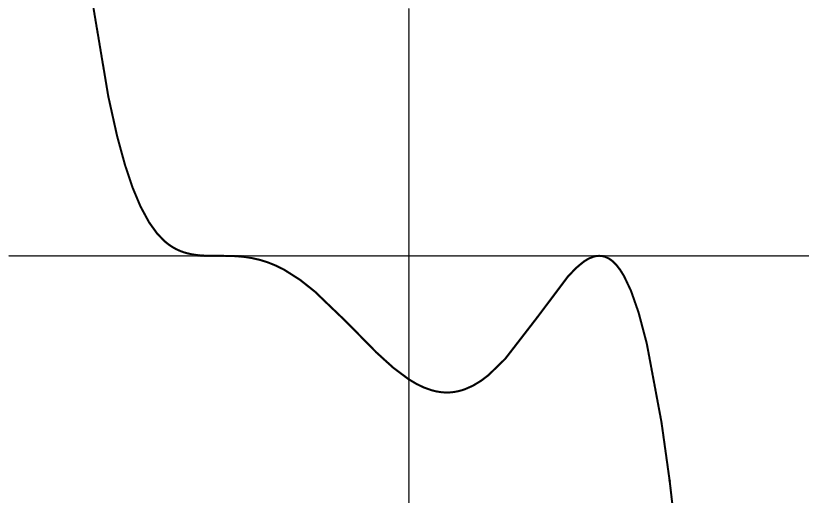}}&
 \resizebox{22mm}{!}{\includegraphics{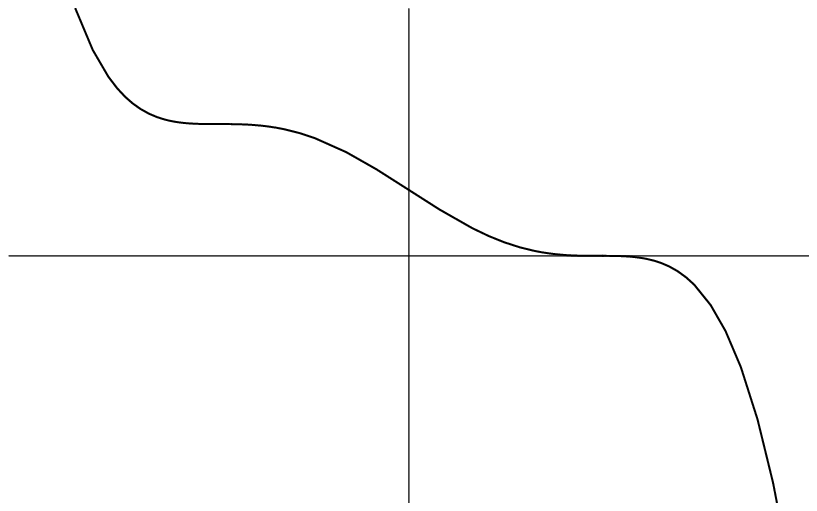}}
  \\
\begin{picture}(1,1.5)
\put(0.5,1.5){\vector(0,-1){1.5}}
\end{picture}&\resizebox{22mm}{!}{\includegraphics{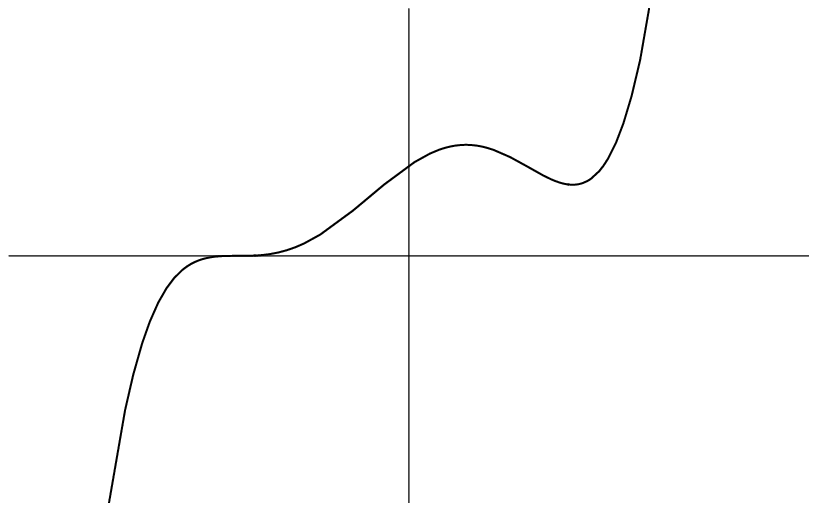}}&
 \resizebox{22mm}{!}{\includegraphics{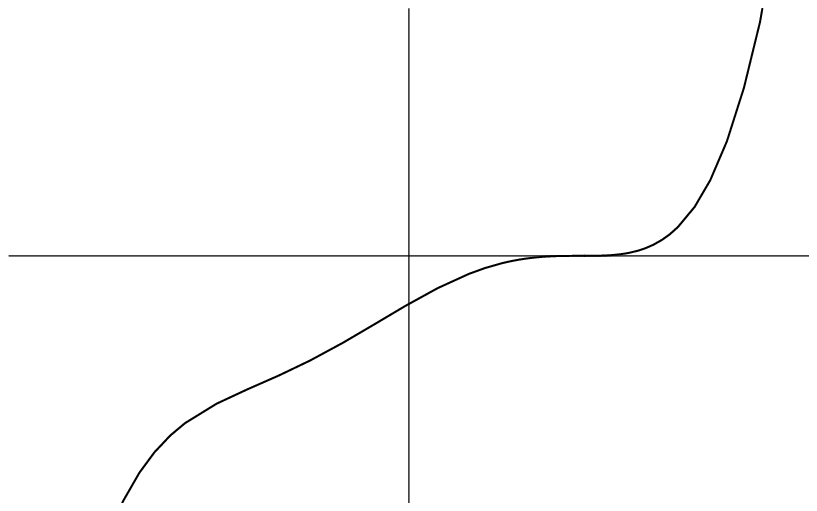}}&
 \resizebox{22mm}{!}{\includegraphics{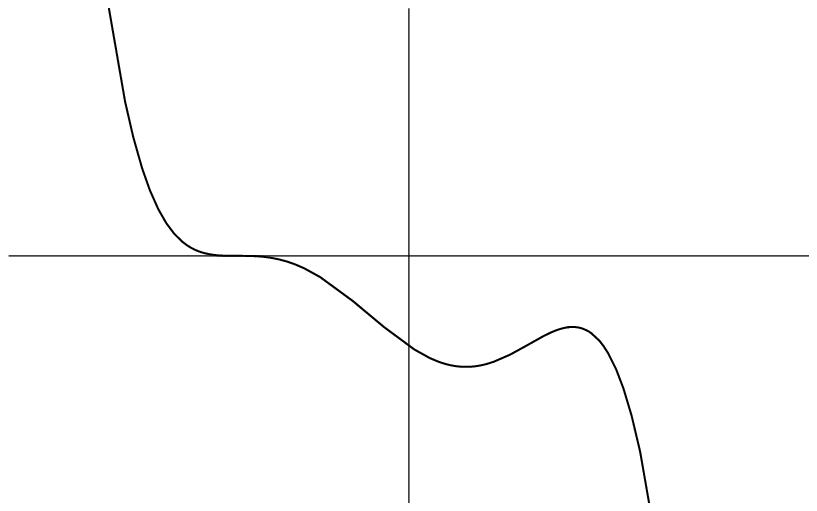}}&
 \resizebox{22mm}{!}{\includegraphics{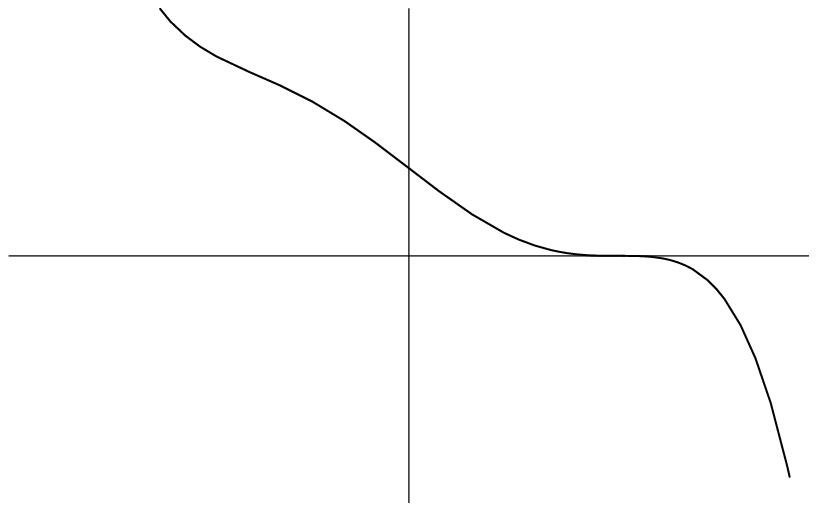}}
  \end{tabular}
 \caption{Graphs of $\mathcal{F}_{\lambda}(x_0^1)$ as $\lambda$ varies.}
 \label{fnmin}
 \end{figure}

\begin{defn}\label{c20}
The pre-normalised reduced action function evaluated on the
caustic is given by
$\mathcal{F}_{\lambda}(x_{0}^{1}): = f_{(x_{t}(\lambda),t)}(x_{0}^{1})-
 f_{(x_{t}(\lambda),t)}(\lambda_{1}).$
 \end{defn}
Assume that $\mathcal{F}_{\lambda}(x_{0}^{1})$ is a real analytic
function in a neighbourhood of $\lambda_{1}\in\mathbb{R}$. Then,
\begin{equation}\label{ce14}\mathcal{F}_{\lambda}(x_{0}^{1}) = (x_{0}^{1}-\lambda_{1})^{3}\tilde{F}(x_{0}^{1}),\end{equation}
where $\tilde{F}$ is real analytic. When the inflexion at $x_0^1=\lambda_1$
 is the minimising
critical point of $\mathcal{F}_{\lambda}$, the caustic will be cool.
Therefore, on a hot/cool boundary this inflexion is about to become
or cease being the minimiser.
\begin{prop}
    A \emph{necessary condition} for
    $x_{t}(\lambda)\in C_{t}$ to be on a hot/cool boundary is that
    either $\tilde{F}(x_{0}^{1})$ or $\tilde{G}(x_{0}^{1})$
        has a repeated root at
    $x_{0}^{1}=r$ where,
    \[\tilde{G}(x_0^1) = 3 \tilde{F}(x_{0}^{1})
    +(x_{0}^{1}-\lambda_{1})\tilde{F}'(x_{0}^{1}).\]
\end{prop}
\begin{proof}
The minimiser could change when either $\tilde{F}$ has a repeated
root which is  the
    minimiser, or there is a second inflexion at a lower minimising value \cite{Neate}.  \end{proof}
The condition is not sufficient as it includes cases where the
minimiser is not about to change (see  Figure \ref{fnmin}).

 \begin{eg}[The polynomial swallowtail]
    Let $V(x,y)\equiv 0$, $k_{t}(x,y)\equiv x$, and
    $S_{0}(x_{0},y_{0}) =
    x_{0}^{5}+x_{0}^{2} y_{0}.$
    This gives global reducibility and $k_{t}(x,y)\equiv x$ means
    that the effect of the noise is to translate $\epsilon=0$ picture
    through $\left(-\epsilon \int_{0}^{t}W_{s}\rmd s,0\right)$.
    A simple calculation gives,
    \[\tilde{F}(x_{0}) = 12\lambda^{2}-3\lambda t+6\lambda x_{0}-tx_{0}+2
    x_{0}^{2},\]
    \[\tilde{G}(x_{0}) = 15\lambda^{2}-4\lambda t+10\lambda x_{0}-2 t
    x_{0}+5x_{0}^{2}.\]
\begin{figure}[h!]
\begin{center}\begin{picture}(8,3.2)
\put(1,-1.2){\resizebox{80mm}{!}{\includegraphics{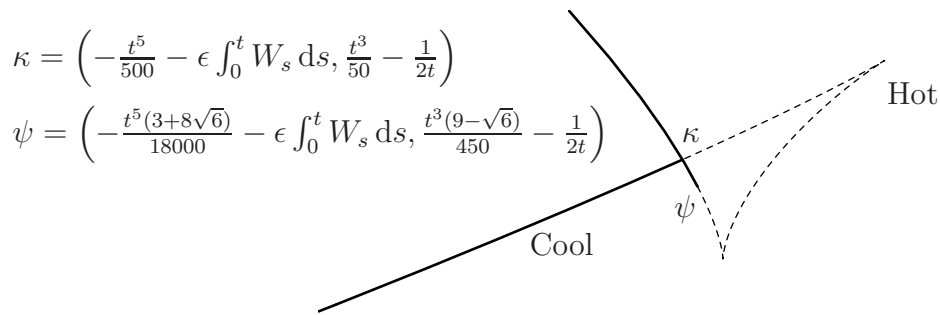}}}
\put(5.8,2){$\kappa $} \put(5.7,1){$\psi$} \put(3.8,0.5){Cool}
\put(-3,3){$\kappa = \left(-\frac{t^5}{500}-\epsilon\int_0^t W_s \rmd s,\frac{t^3}{50}-\frac{1}{2t}\right) $}
\put(-3,2.){$\psi = \left(-\frac{t^5(3+8\sqrt{6})}{18000}-\epsilon\int_0^t W_s\rmd s,\frac{t^3(9-\sqrt{6})}{450}-\frac{1}{2t}\right)$}
\put(8.5,2.5){Hot}
\end{picture}\end{center}
 \caption{Hot and cool parts of the polynomial swallowtail caustic for $t=1$.}\label{swhc}
 \end{figure}
\end{eg}
\begin{eg}[The three dimensional polynomial swallowtail]
 Let $V(x,y)\equiv 0$, $k_{t}(x,y)\equiv 0$, and
    $S_{0}(x_{0},y_{0},z_0) =
    x_{0}^{7}+x_{0}^{3} y_{0}+x_{0}^{2} z_{0}.$
    The functions $\tilde{F}$ and $\tilde{G}$ can be easily found, and an exact expression for the boundary extracted \cite{Neate2} this is  shown in Figure \ref{ps3dhca}.
\begin{figure}[h!]
\begin{tabular}{cc}
\resizebox{50mm}{!}{\includegraphics{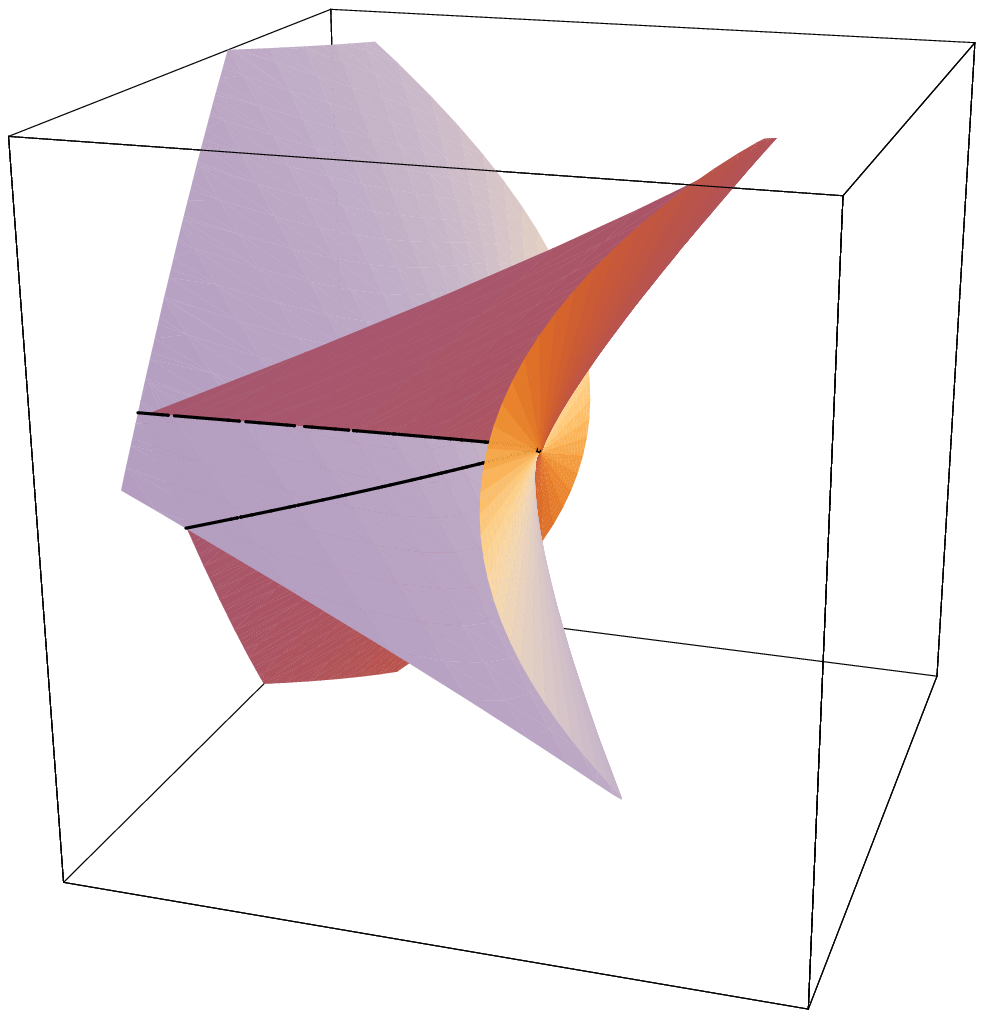}}&
\resizebox{50mm}{!}{\includegraphics{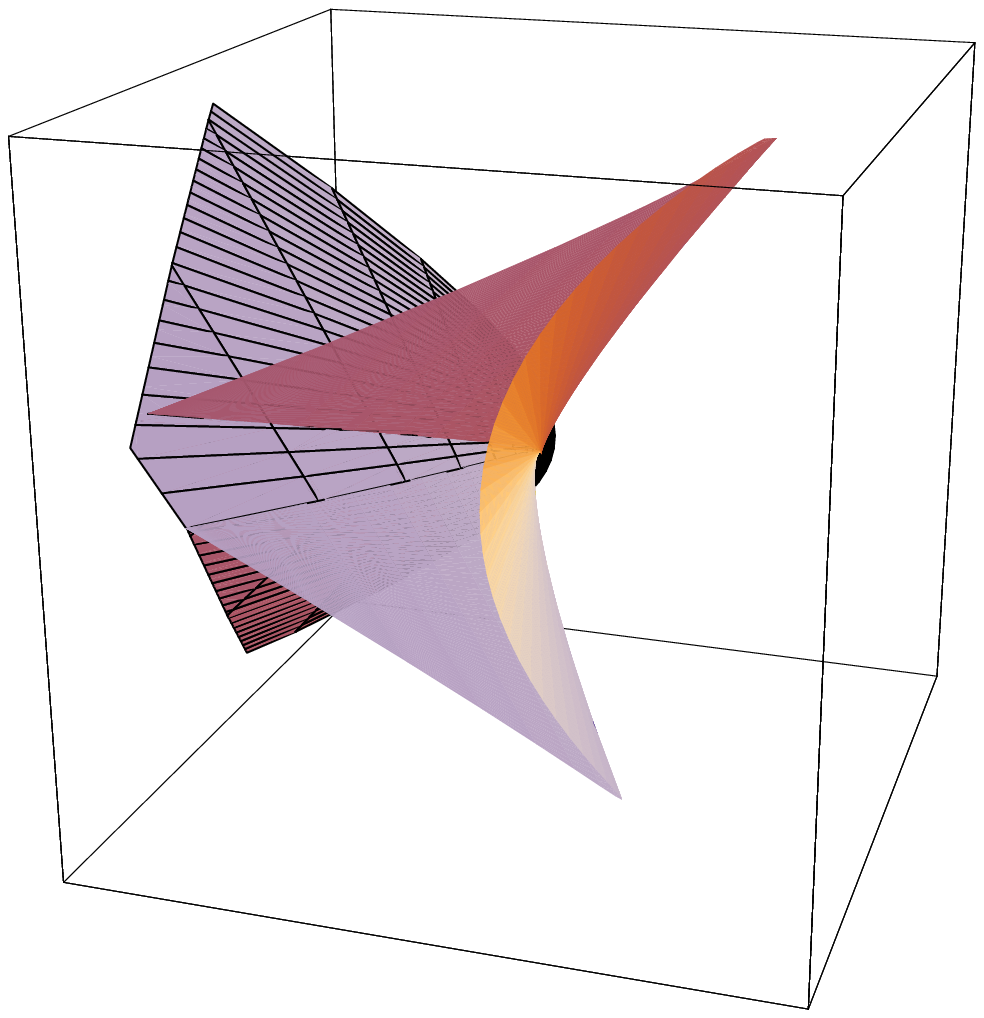}}\\
\emph{Boundary on the caustic. }& \emph{Hot and cool parts.}
\end{tabular}
\caption{The hot (plain) and cool (mesh) parts of the 3D polynomial swallowtail caustic at time $t=1$.}
\label{ps3dhca}
\end{figure}

\end{eg}

\section{Swallowtail perestroikas}
The geometry of a caustic or wavefront can suddenly change with
singularities appearing and disappearing \cite{Arnold3}. We consider
the formation or collapse of a swallowtail using some earlier works
of Cayley and Klein. This section provides a summary of results from
\cite{Neate} where all proofs can be found.

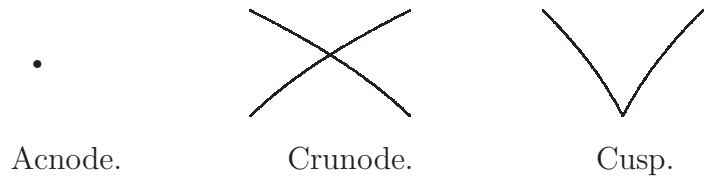
\begin{figure}[h!]
\begin{center}
\setlength{\unitlength}{7mm}
\begin{picture}(15,3)
\put(2,2){\circle*{0.15}}
\qbezier(6,1)(7,2)(9,3) \qbezier(6,3)(8,2)(9,1)
\qbezier(11.5,3)(12.5,2)(13,1) \qbezier(13,1)(13.5,2)(14.5,3)
\put(1.5,0){Acnode.} \put(6.7,0){Crunode.} \put(12.5,0){Cusp.}
\end{picture}
\end{center}
\caption{The classification of double points.} \label{double}
\end{figure}

We begin by recalling the classification of double points of a two
dimensional algebraic curve as acnodes, crunodes and cusps (Figure
\ref{double}).

\begin{figure}[h]
  \setlength{\unitlength}{1 cm} \centering
 \resizebox{110mm}{!}{ \includegraphics{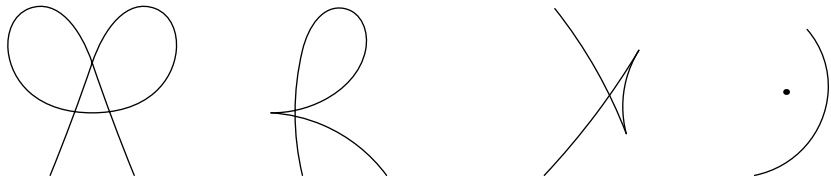}}
  \caption{Cayley's triple points.}
  \label{pic3}
  \end{figure}

In Cayley's work on plane algebraic curves, he describes the
possible triple points of a curve \cite{Salmon} by considering the
collapse of systems of double points which would lead to the
existence of three tangents at a point. The four possibilities are
shown in Figure \ref{pic3}. The systems will collapse to form a
triple point with respectively, three real distinct tangents,
 three real tangents with two coincident,
 three real tangents all of which are coincident,
or one real tangent and two complex tangents. It is the interchange
between the last two cases which will lead to the formation of a
swallowtail on a curve \cite{Hwa}. This interchange was investigated
by Felix Klein \cite{Klein}.

In Section 3, we restricted the pre-parameter to be real to only consider points with real pre-images. This does not allow
there to be any isolated double points. We now allow the parameter to vary
throughout the complex plane and consider when this maps to
real points. We begin by working with a general curve  of the form ${x}(\lambda)
=(x_{1}(\lambda),x_{2}(\lambda))$ where each $x_{\alpha}(\lambda)$ is
real analytic in $\lambda\in\mathbb{C}$. If Im$\left\{{x}(a+\rmi\eta)\right\}=0$,
 it follows that ${x}(a+\rmi\eta)={x}(a-\rmi\eta),$
  so this is a
  ``complex double point'' of the curve ${x}(\lambda)$.

\begin{lemma}
    If ${x}(\lambda) =(x_{1}(\lambda),x_{2}(\lambda))$ is a real
    analytic parameterisation of a curve and $\lambda$ is an
    intrinsic parameter, then there is a generalised cusp at
    $\lambda=\lambda_{0}$ if and only if the curves,
    \[0=\frac{1}{\eta}\mbox{Im} \left\{x_{\alpha}(a+\rmi\eta)\right\}\qquad \alpha=1,2,\]
    intersect at $(\lambda_{0},0)$ in the $(a,\eta)$ plane.
    \end{lemma}

   Now consider a family of parameterised curves
  ${x}_{t}(\lambda)=(x_t^1(\lambda),x_t^2(\lambda))$. As $t$ varies the geometry of the
  curve can change with swallowtails forming and disappearing.

  \begin{prop}\label{swall}
      If a swallowtail on the curve ${x}_{t}(\lambda)$
      collapses to a point
      where $\lambda=\tilde{\lambda}$ when $t=\tilde{t}$ then,
      \[\frac{\rmd {x}_{\tilde{t}}}{\rmd \lambda}(\tilde{\lambda})=
      \frac{\rmd^{2} {x}_{\tilde{t}}}{\rmd \lambda^{2}}(\tilde{\lambda})=0.\]
  \end{prop}

\begin{prop}\label{comp}
Assume that there exists a neighbourhood of $\tilde{\lambda}\in\mathbb{R}$
such that $\frac{\rmd {x}_{t}^{\alpha}}{\rmd \lambda}(\lambda)\neq 0$
for $t\in(\tilde{t}-\delta,\tilde{t})$ where $\delta>0$.
    If a complex double point joins the curve ${x}_{t}(\lambda)$
    at $\lambda=\tilde{\lambda}$ when $t=\tilde{t}$ then,
    \[\frac{\rmd {x}_{\tilde{t}}}{\rmd \lambda}(\tilde{\lambda})=
      \frac{\rmd^{2} {x}_{\tilde{t}}}{\rmd \lambda^{2}}(\tilde{\lambda})=0.\]

\end{prop}

These provide a necessary condition for the
formation or destruction of a swallowtail,  and for complex double
points to join or leave the main curve.
 \begin{defn}
 A  family of parameterised curves $x_t(\lambda)$, (where $\lambda$ is some intrinsic parameter) for which,
 \[\frac{\rmd x_{\tilde{t}}}{\rmd \lambda}(\tilde{\lambda})=
 \frac{\rmd^2 x_{\tilde{t}}}{\rmd \lambda^2}(\tilde{\lambda})=0\]
  is said to have  a point of swallowtail perestroika when $\lambda=\tilde{\lambda}$ and $t=\tilde{t}$.
 \end{defn}
As with generalised cusps,
we have not ruled out further degeneracy at these points.
Moreover, as Cayley highlighted, these points are not cusped and are barely distinguishable from an ordinary point of the curve \cite{Salmon}.

\subsection{The complex caustic in two dimensions}
The complex caustic is the complete caustic found by allowing the parameter $\lambda$ in the  pre-parameterisation  ${x}_{t}(\lambda)\in\mathbb{R}^2$ to vary over the complex plane.
By considering the complex caustic, we are determining
solutions $a=a_{t}$ and $\eta=\eta_{t}$  to,
\[f'_{({x},t)}(a+\rmi\eta) = f''_{({x},t)}(a+\rmi\eta) = 0,\]
where ${x}\in\mathbb{R}^{2}$.
We are interested in these points if they
join the main caustic at some finite critical time $\tilde{t}$. That is,
there exists a finite value $\tilde{t}>0$ such that
$\eta_{t}\rightarrow 0$ as $t\uparrow \tilde{t}$. If this holds then a
swallowtail can develop at the critical time $\tilde{t}$.

\begin{theorem} \label{four}For a two dimensional caustic, assume that
${x}_{t}(\lambda)$ is a real analytic function. If at a time
$\tilde{t}$ a swallowtail perestroika occurs on the caustic, then $x=x_{\tilde{t}}(\lambda)$
is a real solution for ${x}$ to,
\[f'_{({x},\tilde{t})}(\lambda)=f''_{({x},\tilde{t})}(\lambda)=
f'''_{({x},\tilde{t})}(\lambda)=f^{(4)}_{({x},\tilde{t})}(\lambda)=0,
\]
where $\lambda=a_{\tilde{t}}$.
\end{theorem}

\begin{theorem} For a two dimensional caustic, assume that
${x}_{t}(\lambda)$ is a real analytic function. If at a time
$\tilde{t}$  there is a real solution for ${x}$ to,
\[f'_{({x},\tilde{t})}(\lambda)=f''_{({x},\tilde{t})}(\lambda)=
f'''_{({x},\tilde{t})}(\lambda)=f^{(4)}_{({x},\tilde{t})}(\lambda)=0,\]
and  the vectors $\nabla_x f'_{(x,\tilde{t})}(\lambda)$ and $
\nabla_x
    f''_{(x,\tilde{t})}(\lambda)$
    are linearly independent, then $x$ is a point of swallowtail
    perestroika on the caustic.
\end{theorem}

\begin{figure}[h!]
\begin{center}
 \resizebox{110mm}{!}{\includegraphics{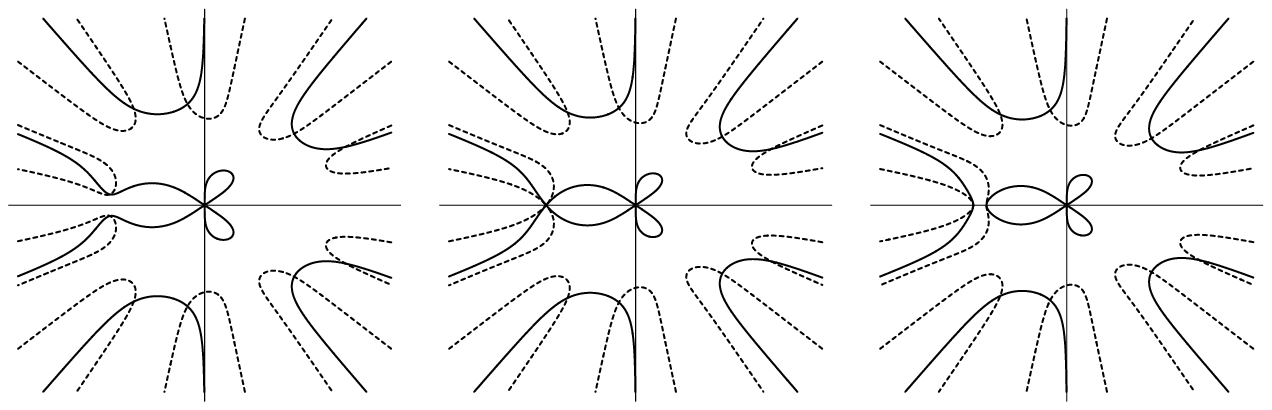}}
 \caption{Im$\left\{x_t(a+\rmi\eta)\right\}=0$ (solid) and Im$\left\{y_t(a+\rmi\eta)\right\}=0$ (dashed) in $(a,\eta)$
plane.} \label{dash1}
\end{center}\end{figure}

\begin{figure}[h!]
 \begin{center}\resizebox{110mm}{!}{\includegraphics{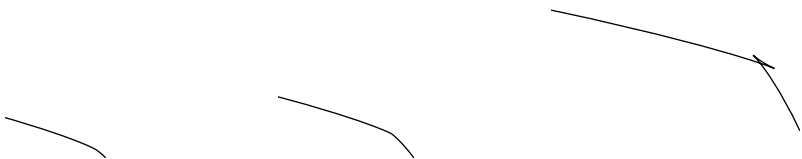}}
\caption{Caustic plotted at corresponding times.}
\label{dash2}
\end{center}
\end{figure}

 \begin{eg}Let $V(x,y)= 0, k_t(x,y)\equiv 0$ and
    $S_{0}(x_{0},y_0) = x_{0}^{5}+x_{0}^{6}y_{0}.$
The caustic  has no cusps for times $t<\tilde{t}$ and two cusps for
times $t>\tilde{t}$ where $\tilde{t} =
    4\sqrt{2}\times 33^{3/4}\times 7^{(-7/4)}
    =2.5854\ldots$

At the critical time $\tilde{t}$ the caustic has a point of swallowtail perestroika
as shown in Figures \ref{dash1} and \ref{dash2}. The conjugate pairs of intersections of the curves in Figure \ref{dash1} are the complex double points. There are five before the critical time and four
afterwards. The remaining complex double points do not join the main
caustic and so do not influence its behaviour for real times.
\end{eg}

\subsection{Level surfaces}
Unsurprisingly, these phenomena are not restricted to caustics. There is an interplay  between the level surfaces
and the caustics, characterised by their pre-images.
\begin{prop}\label{s13}
    Assume that in two dimensions at
    $x_0\in\Phi_t^{-1}H_t^c\cap\Phi_t^{-1}C_t$
    the normal to the pre-level surface $n(x_0)\neq 0$ and the normal
    to the pre-caustic $\tilde{n}(x_0)\neq 0$ so that the pre-caustic
    is not cusped at $x_0$.
    Then  $\tilde{n}(x_0)$ is parallel to $n(x_0)$ if and only if
    there is a generalised cusp on the caustic.
    \end{prop}
    \begin{cor}\label{s14}
    Assume that in two dimensions at
    $x_0\in\Phi_t^{-1}H_t^c\cap\Phi_t^{-1}C_t$
    the normal to the pre-level surface $n(x_0)\neq 0$.
    Then at $\Phi_t(x_0)$ there
    is a point of swallowtail perestroika on the level surface
    $H_t^c$ if and only if  there is a generalised cusp on the caustic $C_t$ at $\Phi_t(x_0)$.
\end{cor}
 \begin{eg}Let $V(x,y)= 0$, $k_t(x,y)=0$, and
 $S_0(x_0,y_0)=x_0^5+x_0^6y_0.$
Consider the behaviour of the level surfaces through a
point inside the caustic swallowtail at a fixed time as the point is moved through a cusp on the
caustic. This is illustrated in Figure  \ref{cuspls}. Part (a) shows
all five of the level surfaces through the point demonstrating how three
swallowtail level surfaces collapse together at the cusp to form a
single level surface with a point of swallowtail perestroika. Parts
(b) and (c) show how one of these swallowtails collapses on its own
and how its pre-image behaves.
\begin{figure}[h!]
\setlength{\unitlength}{10mm}
\begin{picture}(12.5,5.5)
\put(0.5,0){\resizebox{120mm}{!}{\includegraphics{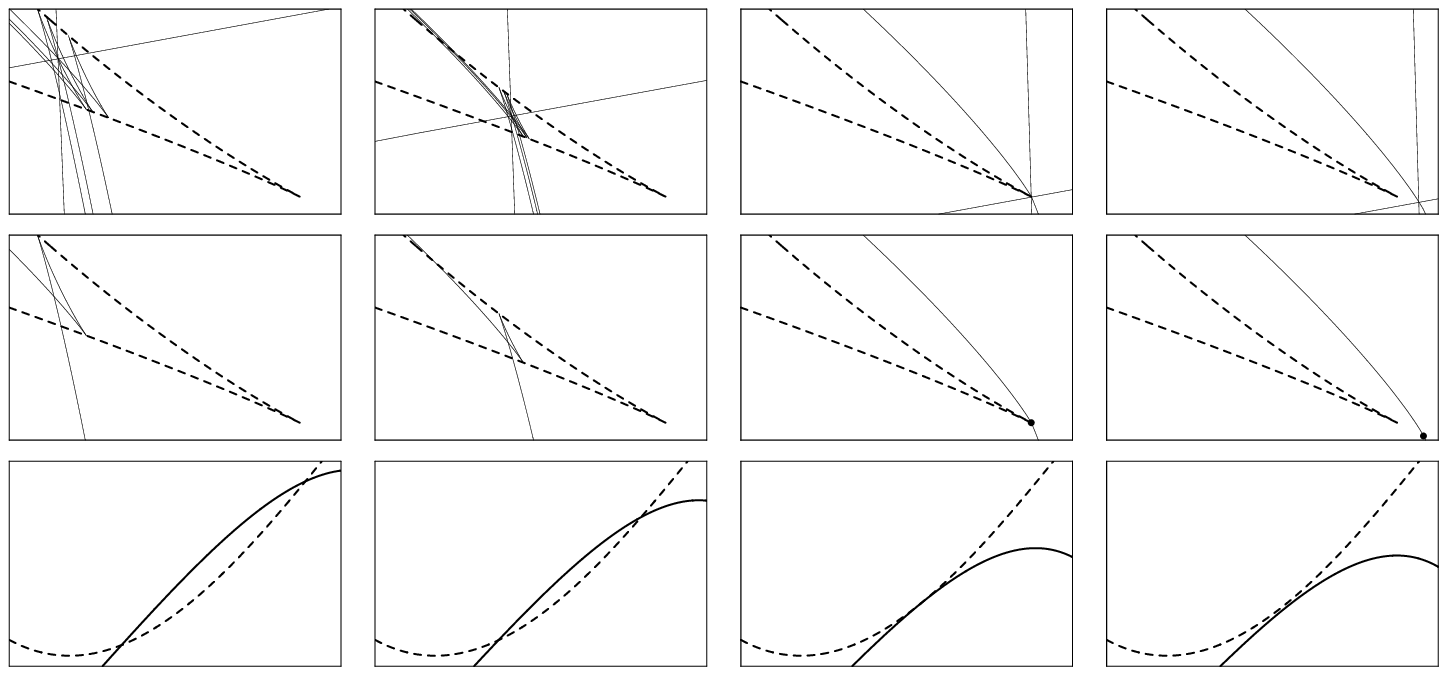}}}
\put(0,1.5){(c)}
\put(0,3.2){(b)}
\put(0,5){(a)}
\end{picture}
\caption{(a) All level surfaces (solid line) through a point as it
crosses the caustic (dashed line) at a cusp, (b) one of these level
surfaces with its complex double point, and (c) its real pre-image.}
\label{cuspls}
 \end{figure}
\end{eg}

\section{Maxwell sets}\label{Mwelsec}

 A jump
will occur in the inviscid limit of the Burgers velocity field if we cross a point at which there are two different
global minimisers $x_0(i)(x,t)$ and $x_0(j)(x,t)$ returning the same value
of the action.

\begin{figure}[h!]
 \setlength{\unitlength}{7.5mm}
\begin{center}
\begin{tabular}{ccc} \textbf{Before Maxwell set}& \textbf{On
Cool Maxwell set} & \textbf{Beyond Maxwell set}\\
\begin{picture}(4.5,3) \put(0,0){\resizebox{33.75mm}{!}{\includegraphics{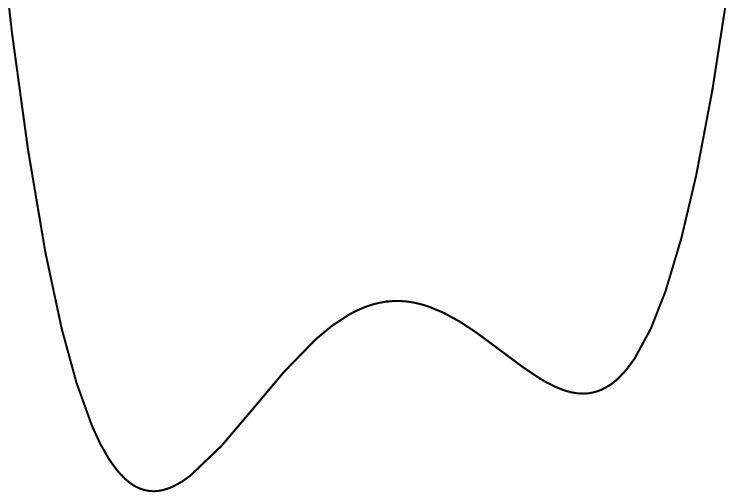}}}
 \put(0.8,0.3){\small{$x_0^1$}}
 \put(3.2,0.9){\small{$\check{x}_0^1$}}\end{picture}&
\begin{picture}(4.5,3) \put(0,0){\resizebox{33.75mm}{!}{\includegraphics{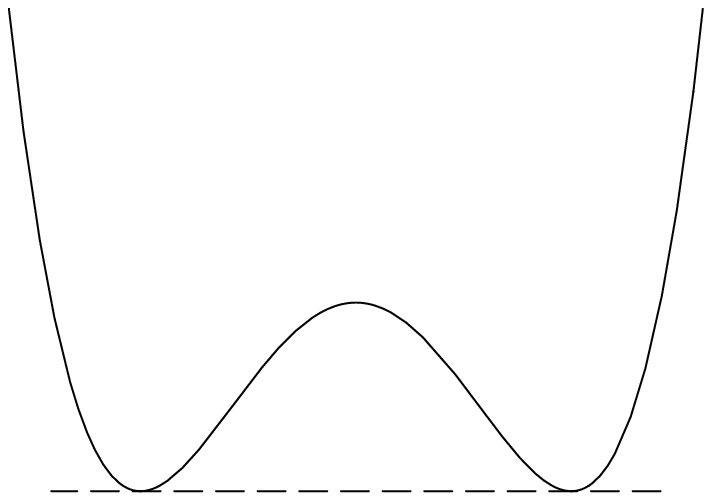}}}
  \put(0.9,0.4){\small{$x_0^1$}}
 \put(3.3,0.4){\small{$\check{x}_0^1$}}\end{picture}&
\begin{picture}(4.5,3)
\put(0,0){\resizebox{33.75mm}{!}{\includegraphics{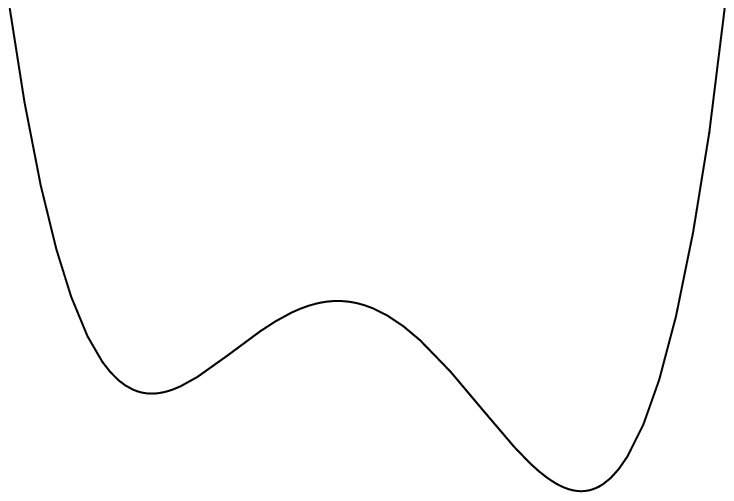}}}
  \put(0.9,0.9){\small{$x_0^1$}}
 \put(3.3,0.3){\small{$\check{x}_0^1$}}
 \end{picture}\\
\emph{ Minimiser at $x_0^1$.}   & \emph{ Two $x_0$'s at same level.}
& \emph{ Minimiser jumps.}
\end{tabular}
\end{center}
\caption{The graph of the reduced action function as $x$ crosses the
Maxwell set.}\label{mf1}
\end{figure}

In terms of the reduced action function, the Maxwell set corresponds
to values of $x$ for which $f_{(x,t)}(x_0^1)$ has two critical
points at the same height. If this occurs at the minimising value
then the Burgers fluid velocity will jump as shown in Figure \ref{mf1}.

%
%
\subsection{The Maxwell-Klein set}\label{ms2}

We begin with the two dimensionals polynomial case by considering
the classification of double points of a curve (Figure
\ref{double}).

\begin{lemma}\label{m4}
A point $x$ is in the Maxwell set if and only if there is a
Hamilton-Jacobi level surface with a point of self-intersection
(crunode) at $x$.
\end{lemma}
\begin{proof} Follows from Definition \ref{m1}.
\end{proof}
\begin{defn}\label{m5}
The Maxwell-Klein set $B_t$ is the set of points which are non-cusp
double points of some Hamilton-Jacobi level surface curve.
\end{defn}

It follows from this definition that a point is in the Maxwell-Klein
set if it is either a complex double point (acnode) or point of self-intersection (crunode) of some Hamilton-Jacobi level surface. Using the geometric results of DTZ outlined in
Section 2, it is easy to calculate this set in the polynomial
case as the cusps of the level surfaces sweep out the caustic.

\begin{theorem}\label{m6}
Let $D_t$ be the set of double points of the
Hamilton-Jacobi level surfaces, $C_t$ the caustic set, and $B_t$ the
Maxwell-Klein set. Then, from Cayley and Klein's classification of
double points as crunodes, acnodes, and cusps, by definition
$D_t=C_t\cup B_t$ and the corresponding defining algebraic equations
factorise $D_t=C_t^n\cdot B_t^m$, where $m,n$ are positive integers.
\end{theorem}
\begin{proof}
Follows from Proposition \ref{i11} and Lemma \ref{m4}.
\end{proof}

\begin{theorem}\label{2dbig2}
Let $\rho_{(t,c)}(x)$ be the resultant,
\[ \rho_{(t,c)}(x)=R\left(f_{(x,t)}(\cdot)-c, f_{(x,t)}'(\cdot)\right),\] where $x=(x_1,x_2)$.
Then $x\in D_t$  if and only if for some $c$,
\[\rho_{(t,c)}(x)=\frac{\partial \rho_{(t,c)}}{\partial x_1}(x)=
\frac{\partial \rho_{(t,c)}}{\partial x_2}(x)=0.\] Further,
\[D_t(x) = \gcd\left(\rho_t^1(x),\rho^2_t(x)\right),\]
 where $\gcd(
\cdot,\cdot )$ denotes the greatest common divisor and $\rho_t^1$ and $\rho_t^2$ are the resultants,
\[{\rho}_t^1(x)=R\left(\rho_{(t,\cdot)}(x),\frac{\partial \rho_{(t,\cdot)}}{\partial x_1}(x)\right)
\quad\mbox{and}\quad
{\rho}^2_t(x)=R\left(\frac{\partial \rho_{(t,\cdot)}}{\partial x_1}(x),\frac{\partial \rho_{(t,\cdot)}}{\partial x_2}(x)\right).\]
\end{theorem}
\begin{proof} Recall that the equation of the level surface of
Hamilton-Jacobi functions is merely the result of eliminating
$x_0^1$ between the equations,
\[f_{(x,t)}(x_0^1)=c\quad\mbox{and}\quad f_{(x,t)}'(x_0^1)=0.\]
We form the resultant $\rho_{(t,c)}(x)$ using Sylvester's formula.  The double
points of the level surface must satisfy for some $c\in\mathbb{R}$,
\[\rho_{(t,c)}(x)=0,\quad\frac{\partial \rho_{(t,c)}}{\partial x_1}(x)=0\quad\mbox{and}\quad
\frac{\partial \rho_{(t,c)}}{\partial x_2}(x)=0.\] Sylvester's
formula proves all three equations are polynomial in $c$. To proceed we
eliminate $c$ between pairs of these equations using resultants
giving,
\[R\left(\rho_{(t,\cdot)}(x),\frac{\partial \rho_{(t,\cdot)}}{\partial x_1}(x)\right)={\rho}_t^1(x)
\quad\mbox{and}\quad R\left(\frac{\partial
\rho_{(t,\cdot)}}{\partial x_1}(x),\frac{\partial
\rho_{(t,\cdot)}}{\partial x_2}(x)\right)={\rho}^2_t(x).\] Let
$D_t=\gcd(\rho_t^1,\rho^2_t)$ be the greatest common divisor of the
algebraic $\rho_t^1$ and $\rho_t^2$. Then $D_t(x)=0$ is the equation of double
points.\end{proof}

We now extend this to $d$-dimensions, where the Maxwell-Klein set corresponds to points which satisfy the Maxwell set condition but have both real  pre-images (Maxwell) or complex  pre-images (Klein).
\begin{theorem}\label{m10}
Let the reduced action function $f_{(x,t)}(x_0^1)$ be a polynomial
in all space variables. Then the set of all possible discontinuities
for a $d$-dimensional Burgers fluid velocity field in the inviscid
limit is the double discriminant,
 $$D(t):=D_c\left\{D_{\lambda}\left(f_{(x,t)}(\lambda)-c\right)\right\}=0,$$
 where $D_x(p(x))$ is the discriminant of the polynomial $p$ with
 respect to $x$.
\end{theorem}
\begin{proof}
By considering the Sylvester matrix of the first discriminant,
\[D_{\lambda}\left(f_{(x,t)}(\lambda)-c\right)=K\prod\limits_{i=1}^m \left(f_{(x,t)}(x^1_0(i)(x,t))-c\right),\]
where $x_0^1(i)(x,t)$ is an enumeration of the real and complex
roots $\lambda$ of $f'_{(x,t)}(\lambda)=0$ and $K$ is some constant.
Then  the second discriminant is simply, \[D_c\left(D_{\lambda}\left(f_{(x,t)}(\lambda)-c\right)\right) =
K^{2m-2}\prod\limits_{i<j}\left(f_{(x,t)}(x^1_0(i)(x,t))-f_{(x,t)}(x^1_0(j)(x,t))\right)^2.\qedhere\]
\end{proof}

\begin{figure}[h!]
\setlength{\unitlength}{0.9cm}
\begin{center}\begin{picture}(8,8)
  \put(0,0){ \resizebox{6.75cm}{!}{\includegraphics{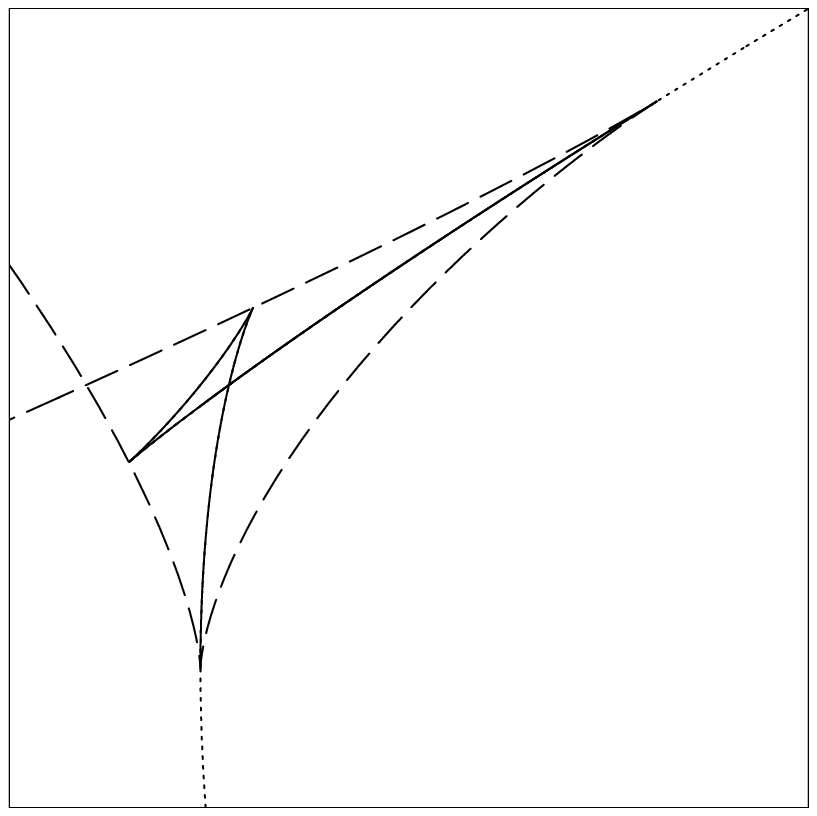}}}
  \put(3,0.8){Acnodes of $H_t^c$}
  \put(4,0.3){ = Klein set.}
  \put(3,2){Crunodes of $H_t^c$}
  \put(4,1.5){ = Maxwell set.}
  \put(0.5,6){Cusps of $H_t^c$ = Caustic.}
  \put(3,0.9){\vector(-1,0){1}}
    \put(2,5.9){\vector(0,-1){1.4}}
      \put(3,2.1){\vector(-1,0){1}}
   \end{picture}\end{center}
\caption{The caustic and Maxwell-Klein set.}\label{mf2}
\end{figure}

\begin{theorem}\label{m13}
The double discriminant $D(t)$ factorises as,
\[D(t) = b_0^{2m-2}\cdot\left(C_t\right)^3\cdot \left(B_t\right)^2,\]
where $B_t=0$ is the equation of the Maxwell-Klein set and $C_t=0$
is the equation of the caustic. The expressions $B_t$ and $C_t$
are both algebraic in
$x$ and $t$.
\end{theorem}
\begin{proof}
See \cite{Neate}.
\end{proof}
 \begin{eg}[The polynomial swallowtail]\label{m9}
 Let $V(x,y)=0$, $k_t(x,y)=0$ and,
 $S_0(x_0,y_0)=x_0^5+x_0^2y_0.$
The Maxwell-Klein set can be found by factorisation giving,
 \begin{eqnarray*}
 0 & = & -675 +52t^4-t^8+3120t^3x-224t^7x +4t^{11}x -38400t^2x^2+1408t^6 x^2 \\
 &   &\quad+128000tx^3
-5400ty+312t^5y-4t^9y+12480t^4xy -448t^8xy \\
&   &\quad-76800t^3x^2y -16200t^2y^2
+624t^6y^2 -4t^{10}y^2 +12480t^5xy^2\\
&   &\quad -21600t^3y^3 +416t^7y^3
-10800t^4y^4.
 \end{eqnarray*}
Outside of the swallowtail on the caustic there are two real and two
complex pre-images whereas inside there are four real and no complex pre-images. Therefore, any part of the Maxwell-Klein set outside of the caustic swallowtail must correspond to Klein double points and any part inside must correspond to the Maxwell set.
This is shown in Figure \ref{mf2}.

 \end{eg}

\subsection{The pre-Maxwell set}\label{ms4}

If the Maxwell set  is defined as in Definition \ref{m1}, then the
pre-Maxwell set is the set of all the pre-images $x_0$ and
$\check{x}_0$ which give rise to the Maxwell set.

 \begin{defn}\label{m15}
The pre-Maxwell set $\Phi_t^{-1}M_t$  is the set of all points
$x_0\in\mathbb{R}^d$ where there exists
$x,\check{x}_0\in\mathbb{R}^d$ such that  $x=\Phi_t(x_0)$ and
$x=\Phi_t(\check{x}_0)$ with $x_0\neq \check{x}_0$ and,
 \[\mathcal{A}(x_0,x,t)=\mathcal{A}(\check{x}_0,x,t).\]
 \end{defn}

With the caustic and level surfaces, each regular point was linked
by $\Phi_t^{-1}$ to a single point on the relevant pre-surafce.
However, every point on the Maxwell set is linked by $\Phi_t^{-1}$
to at least two points on the pre-Maxwell set.
\begin{theorem}\label{m17}
The pre-Maxwell set is given by the discriminant
$D_{\check{x}_0^1}\left(G(\check{x}_0^1)\right)=0$
where, \[G(\check{x}_0^1)=\frac{f_{(\Phi_t (x_0),t)}(x_0^1)-f_{(\Phi_t (x_0),t)}(\check{x}_0^1)}{(x_0^1-\check{x}_0^1)^2}.\]
\end{theorem}
\begin{proof}
From the Definition \ref{m15} and Theorem \ref{i19} it follows that
the pre-Maxwell set is found by eliminating $x$ and $\check{x}_0^1$
between, \[ f_{(x,t)}(x_0^1)=f_{(x,t)}(\check{x}_0^1) \quad
f_{(x,t)}'(x_0^1)=f_{(x,t)}'(\check{x}_0^1)=0.\]This surface would
include the pre-caustic where $x_0^1 = \check{x}_0^1$ and so this
repeated root must be eliminated.
\end{proof}

We can use this to pre-parameterise the Maxwell set as has been done
with the caustic and level surfaces. By restricting the parameter to
be real, we only get the Maxwell set as the Klein points have
complex pre-images.

We now summarise the results of \cite{Neate3}.

\begin{lemma}\label{m22}
Assume that a point $x$ on the Maxwell set corresponds to exactly two pre-images on the pre-Maxwell set, $x_0$ and $\tilde{x}_0$.
Then the normal to the pre-Maxwell set at $x_0$ is to within a scalar multiplier given by,
\begin{eqnarray*}
n(x_0)&  = &-
\left(\frac{\partial^2\mathcal{A}}{\partial x_0^2}(x_0,x,t)\right)
\left(\frac{\partial^2\mathcal{A}}{\partial x\partial x_0}(x_0,x,t)\right)^{-1}
\\
&  &\qquad\qquad\left(\dot{X}(t,x_0,\nabla S_0(x_0))-\dot{X}(t,\tilde{x}_0,\nabla S_0(\tilde{x}_0))\right).\end{eqnarray*}
\end{lemma}

 \begin{cor}\label{m23}
 In two dimensions let the pre-Maxwell set meet the pre-caustic at a point $x_0$ where $n\neq 0$ and
 \[\ker\left(\frac{\partial^2\mathcal{A}}{(\partial x_0)^2}(x_0,\Phi_t(x_0),t)\right)=\langle e_0\rangle,\]
 where $e_0$ is the zero eigenvector. Then the tangent plane to the pre-Maxwell set at $x_0$, $T_{x_0}$ is spanned by $e_0$.
 \end{cor}
\begin{prop}\label{m24}
Assume that in two dimensions at $x_0\in\Phi_t^{-1}M_t$ the normal $n(x_0)\neq 0$ so that the pre-Maxwell set does not have a generalised cusp at $x_0$. Then the Maxwell set can only have a cusp at $\Phi_t(x_0)$ if $\Phi_t(x_0)\in C_t$. Moreover, if \[x=\Phi_t(x_0)\in\Phi_t\left\{ \Phi_t^{-1}C_t\cap \Phi_t^{-1}M_t\right\},\] the Maxwell set will have a generalised cusp at $x$.
\end{prop}

\begin{cor}\label{m25}
In two dimensions, if the pre-Maxwell set intersects the pre-caustic at a point $x_0$, so that there is a cusp on the Maxwell set at the corresponding point where it intersects the caustic,  then the pre-Maxwell set touches the pre-level surface $\Phi_t^{-1} H_t^c$ at the point $x_0$. Moreover, if the cusp on the Maxwell set intersects the caustic at a regular point of the caustic, then there will be a cusp on the pre-Maxwell set which also meets the same pre-level surface $\Phi_t^{-1} H_t^c$ at another point $\check{x}_0$.\end{cor}
\begin{cor}\label{m26}
When the   pre-Maxwell set touches the pre-caustic and pre-level
surface, the Maxwell set intersects a cusp on the caustic.
\end{cor}

\begin{eg}[The polynomial swallowtail]
Let $V(x,y)=0$, $k_t(x,y)=0$ and,
 $S_0(x_0,y_0)=x_0^5+x_0^2y_0.$


\begin{figure}[h]
\setlength{\unitlength}{8mm}
\centering \begin{tabular}{cc}
 \begin{picture}(7,7)
 \put(2.2,4.6){1}
  \put(0.4,5.2){2}
 \put(1.7,5.8){3}
  \put(4.8,5.3){4}
  \put(6.4,0.2){5}
 \put(4.8,3){6}
\put(0,0){\resizebox{56mm}{!}{\includegraphics{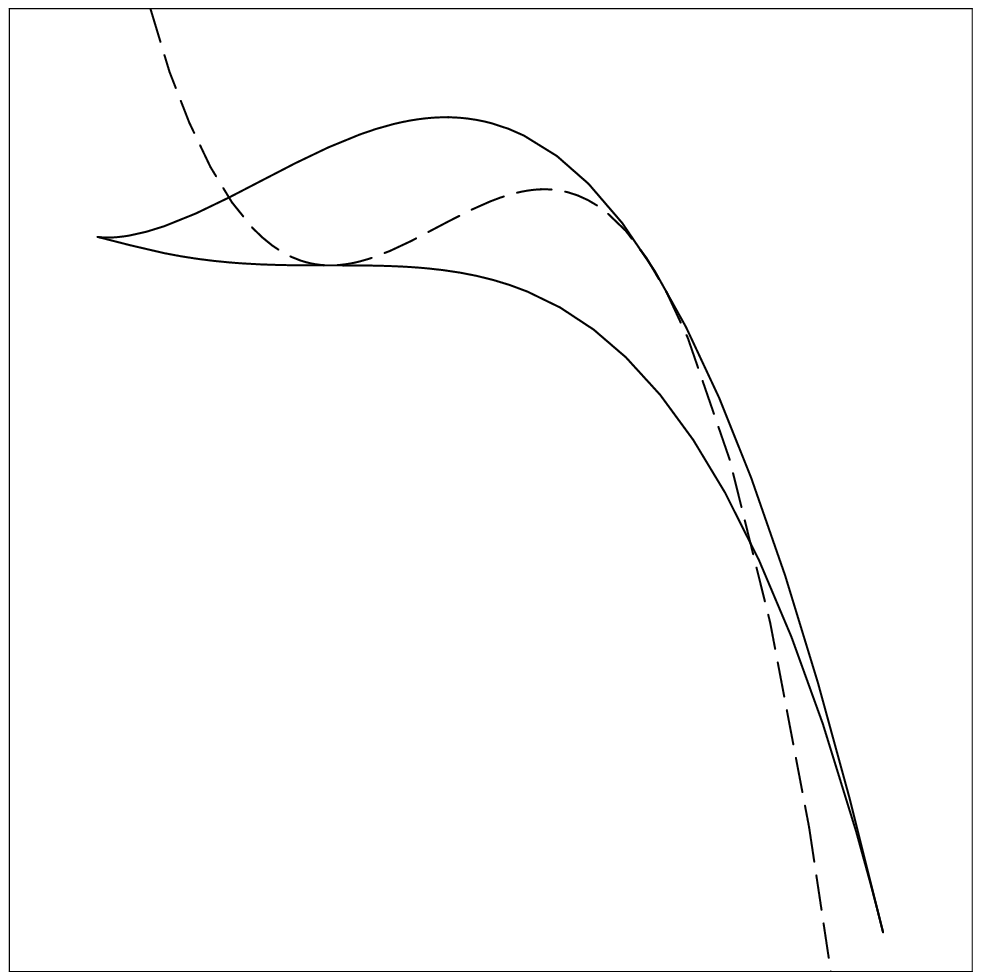}}}
\end{picture}
&
 \begin{picture}(7,7)
 \put(2,1){1}
  \put(2.4,4.2){2}
 \put(1.2,2.8){3}
  \put(6.2,6.0){4}
  \put(1.4,2.5){5}
 \put(2.8,4.3){6}
\put(0,0){\resizebox{56mm}{!}{\includegraphics{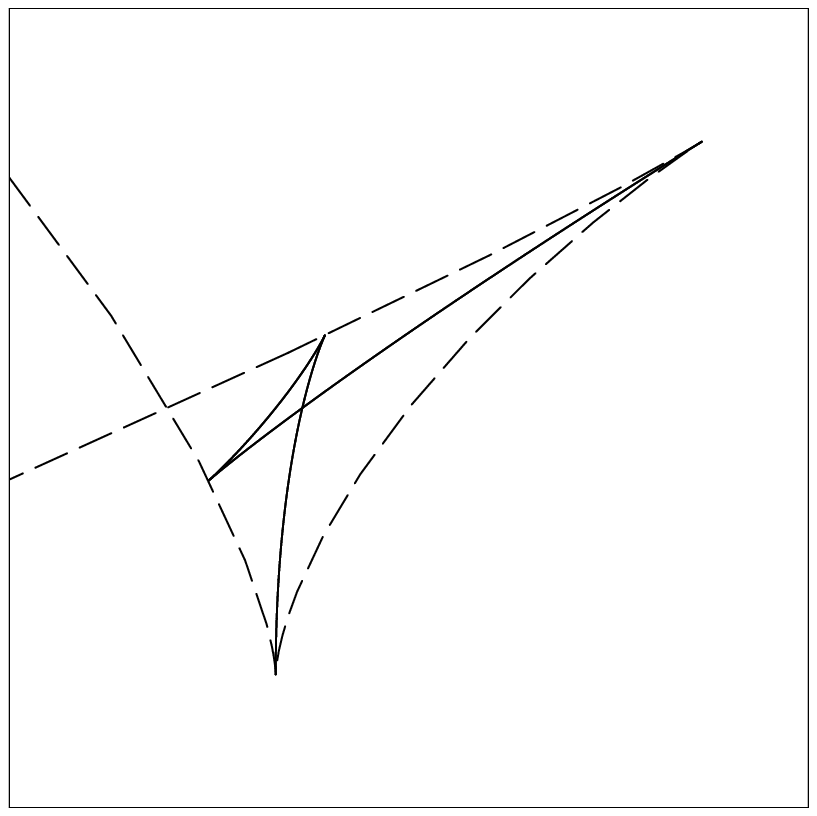}}}
\end{picture}
\\
\emph{ Pre-curves}   & \emph{ Curves } \\
\end{tabular}
\caption{The caustic (dashed) and Maxwell set (solid line).} \label{mf4}
\end{figure}
From Proposition \ref{m24}, the cusps on the Maxwell set correspond to the intersections of the pre-curves (points 3 and 6 on Figure \ref{mf4}). But from Corollary \ref{m25}, the cusps on the Maxwell set also correspond to the cusps on the pre-Maxwell set (points 2 and 5 on Figure \ref{mf4} and also Figure \ref{mf4a}). The Maxwell set terminates when it reaches the cusps on the caustic. These points satisfy the condition for a generalised cusp but, instead of appearing cusped, the curve stops and the parameterisation begins again in the sense that it maps back exactly on itself. At such points the pre-surfaces all touch (Figure \ref{mf4a}).

\begin{figure}[h]
\centering \begin{tabular}{cc}
\resizebox{58mm}{!}{\includegraphics{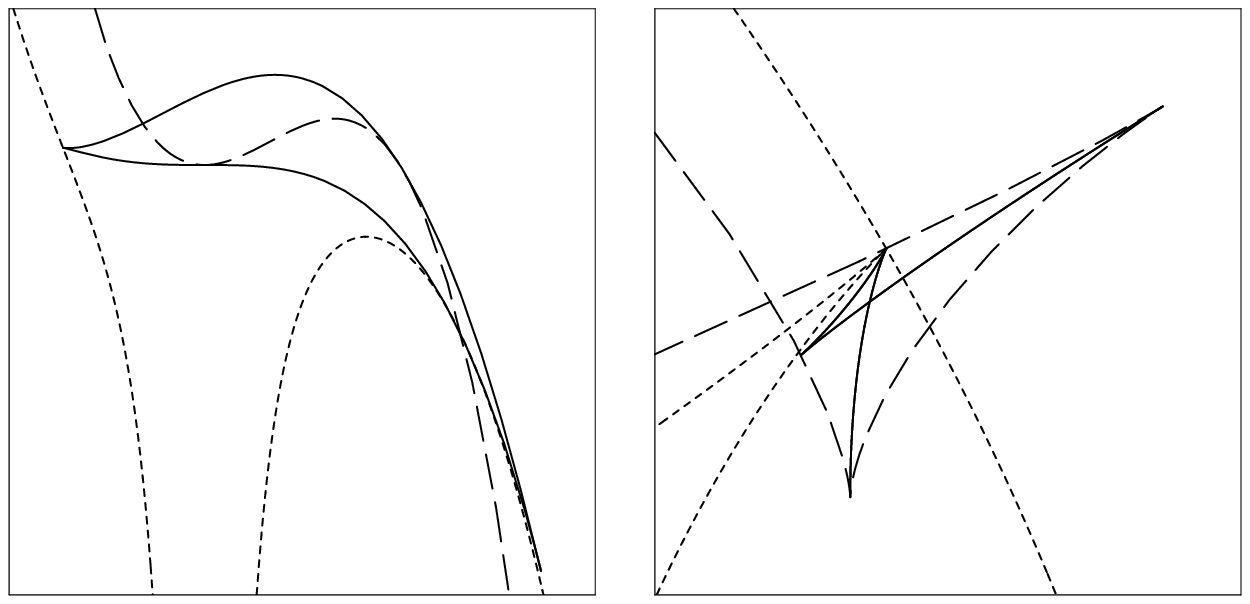}}&
\resizebox{58mm}{!}{\includegraphics{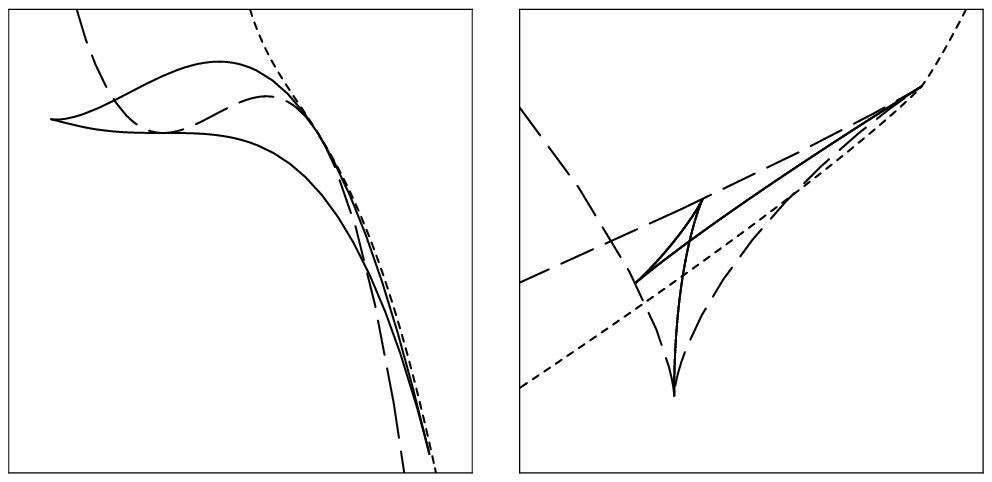}}\\
\emph{ Cusp on Maxwell set}   & \emph{ Cusp on caustic } \\
\end{tabular}
\caption{The caustic (long dash) and Maxwell set (solid line) with the level surfaces (short dash) through special points.} \label{mf4a}
\end{figure}

These two different forms of cusps correspond to very different geometric behaviours of the level surfaces. Where the Maxwell set stops or cusps corresponds to the disappearance of a point of self-intersection on a level surface.
There are two distinct ways in which this can happen. Firstly, the level surface will have a point of swallowtail perestroika when it meets a cusp on the caustic. At such a point only one point of self-intersection will disappear, and so there will be only one path of the Maxwell set which will terminate at that point. However, when we approach the caustic at a regular point, the level surface  must have a cusp but not a  swallowtail perestoika. This corresponds to the collapse of the second system of double points in Figure
\ref{pic3}. Thus, two different points of self intersection coalesce and
so two paths
of the Maxwell set must approach
the point and produce the cusp
(see Figure \ref{mf5}).

  \begin{figure}[h!]\centering
\begin{tabular}{cccc}
\resizebox{27mm}{!}{\includegraphics{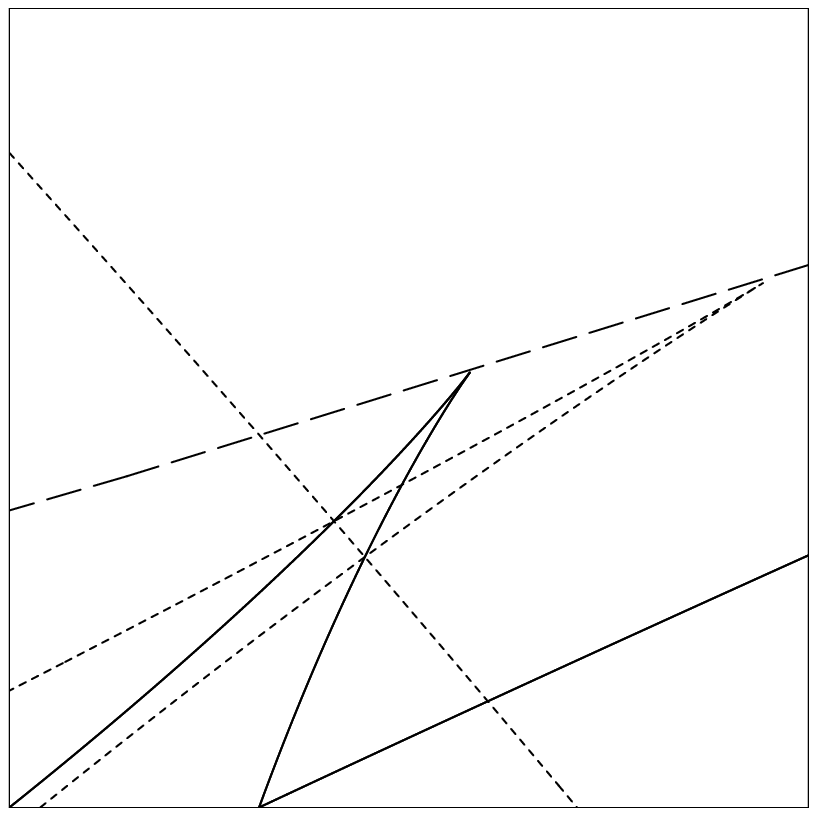}}&
\resizebox{27mm}{!}{\includegraphics{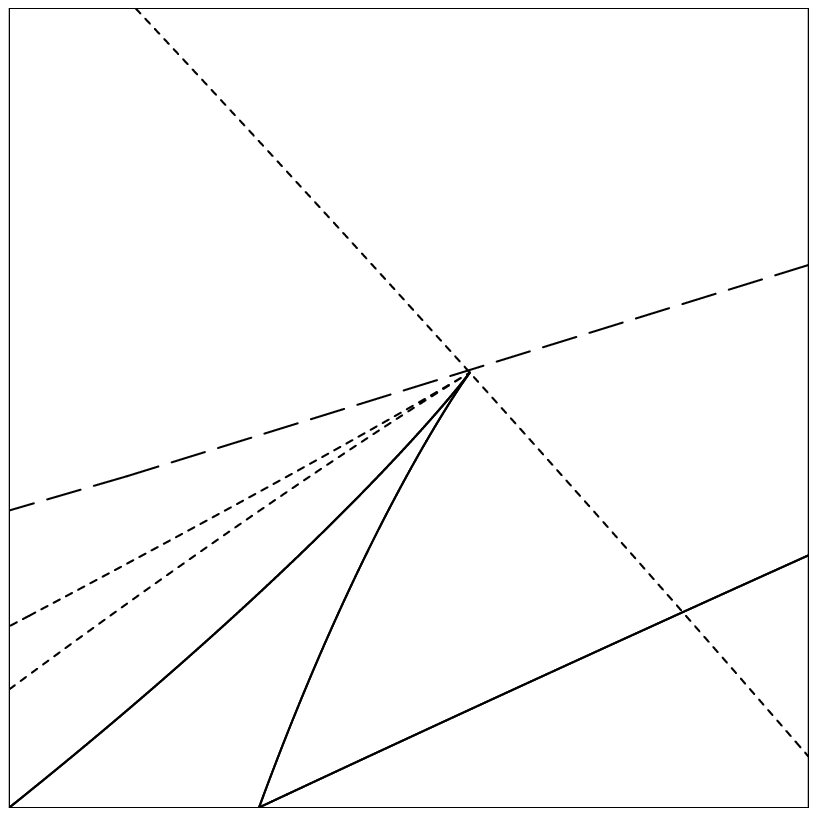}}&
\resizebox{27mm}{!}{\includegraphics{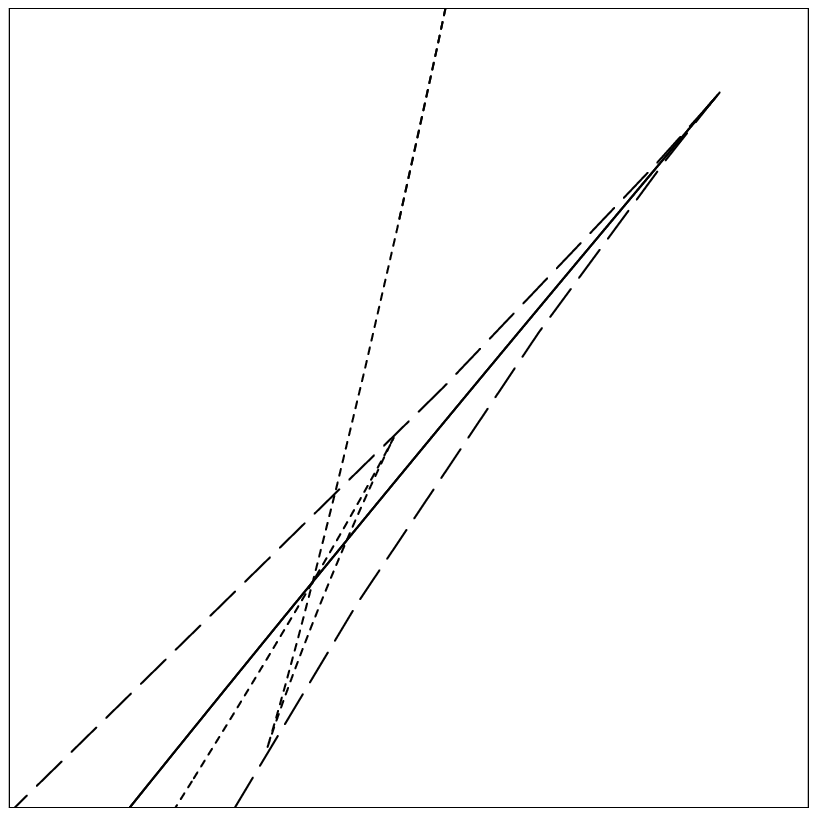}}&
\resizebox{27mm}{!}{\includegraphics{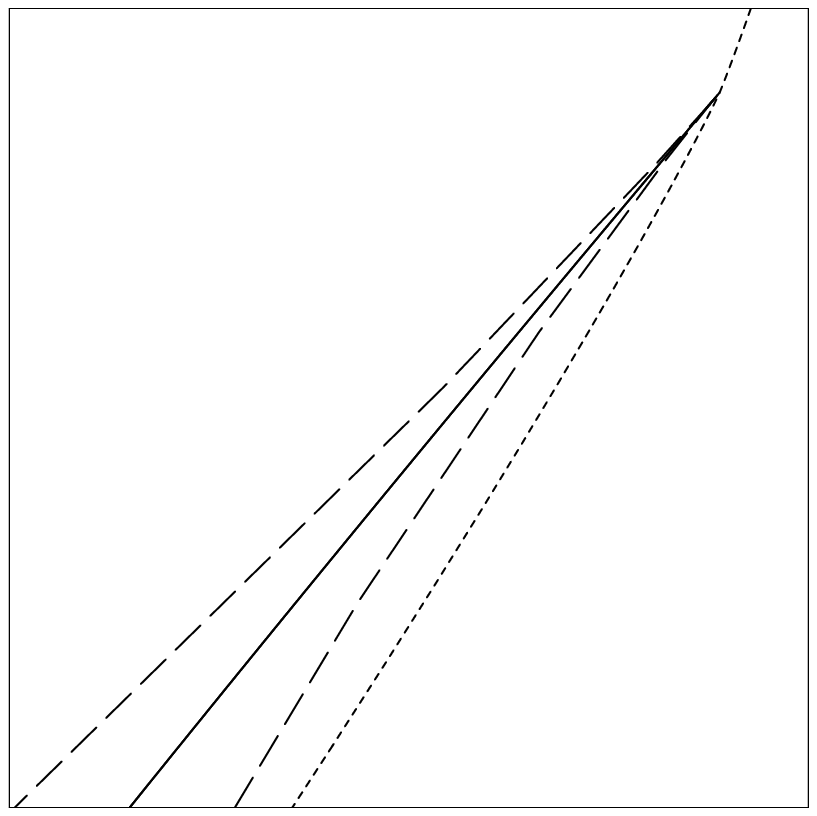}}
\\
\multicolumn{2}{c}{\emph{ Approaching the caustic } }  & \multicolumn{2}{c}{\emph{ Approaching a cusp on the caustic }}
\end{tabular}
\caption{The caustic (long dash) Maxwell set (solid line) and level surface (short dash).} \label{mf5}
\end{figure}
\end{eg}

\section{Some applications to turbulence in two dimensions}

\subsection{Real turbulence and the $\zeta$ process }
\begin{defn}
    The turbulent times $t$ are times when the pre-level surface of the
    minimising Hamilton-Jacobi function \emph{touches}
    the pre-caustic. Such
    times $t$ are zeros of a stochastic process $\zeta^{c}(.)$. i.e. $\zeta^c(t)=0$.
\end{defn}
 These turbulent times are times at which the number of cusps on the corresponding
 level surface will change.
 We begin with some minor generalisations of results in RTW \cite{Truman6} and also \cite{Neate,Reynolds}.
\begin{prop}
    Assume $\Phi_{t}$ is globally reducible and that $x_t(\lambda)$ is the pre-para\-meter\-isation of a two dimensional caustic.
    Then the turbulence process at $\lambda$ is given by,
      \[\zeta^{c}(t) = f_{(x_{t}(\lambda_{0}),t)}(\lambda_{0})-c,\]
    where $f_{(x,t)}(x_{0}^{1})$ is the
    \emph{reduced action} evaluated at
    points
    $x=x_{t}(\lambda_{0})$ where $x_{t}(\lambda_{0})=
    \Phi_{t}(\lambda_{0},x_{0}^{2}(\lambda_{0}))
    \in C_{t},$ $\lambda = \lambda_{0}$ satisfying,
     \[\dot{X}_{t}(\lambda)\cdot\frac{\rmd x_{t}}{\rmd \lambda}(\lambda)
    =0,\]
    where $\dot{X}_t(\lambda) =
    \dot{\Phi}_{t}(\lambda,x_{0,\mathrm{C}}^{2}(\lambda))$ and
    $x_{t}(\lambda_{0})\in C_{t}^{\mathrm{c}}$,
    the cool part of the caustic.
\end{prop}
Hence, there are three kinds of real stochastic turbulence:-
\begin{enumerate}
    \item  \emph{Cusped}, where there is a cusp on the caustic,
    \item \emph{Zero speed}, where the Burgers fluid velocity is zero,
    \item \emph{Orthogonal}, where the Burgers fluid velocity is orthogonal
    to the caustic.
\end{enumerate}
\begin{proof}
The number of cusps on the relevant pre-level surface is,
\[n_{c}(t) =
\#\left\{\lambda\in\mathbb{R}:f_{(x_{t}(\lambda),t)}(\lambda)
=c\right\},\] where the roots $\lambda=\lambda_{0}$ correspond to
points in the cool part of the caustic. The pre-surfaces touch when
$n_{c}(t)$ changes, which occurs when,
\[\frac{\rmd}{\rmd \lambda}f_{(x_{t}(\lambda),t)}(\lambda)=0.\qedhere\]
\end{proof}

For stochastic turbulence to be intermittent we require that
the process
$\zeta^{c}(t)$ is recurrent.
\begin{prop}
    Let $V(x,y)=0$, $k_{t}(x,y) = x$ and \[S_{0}(x_{0},y_{0}) =
    f(x_{0})+g(x_{0}) y_{0},\] where $f,g,f'$ and $g'$ are zero at
    $x_{0}=a$ but $g''(a)\neq 0$. Then, for orthogonal turbulence at $a$,
    \[\zeta^{c}(t) = -a\epsilon W_{t} +\epsilon^{2}
    W_{t}\int_{0}^{t}W_{s}\rmd s -\hf[\epsilon^{2}]
    \int_{0}^{t}W_{s}^{2}\rmd s -c.\]
\end{prop}
We note the following result of RTW \cite{Truman6}.
\begin{lemma}
    Let $W_t$ be a $BM(\mathbb{R})$ process starting at $0$, $c$ any
    real constant and
    \[Y_{t} =  -a\epsilon W_{t} +\epsilon^{2}
    W_{t}\int_{0}^{t}W_{s}\rmd s -\hf[\epsilon^{2}]
    \int_{0}^{t}W_{s}^{2}\rmd s -c.\]
    Then, with probability one, there exists a sequence of times
    $t_{n}\nearrow \infty$ such that
    \[Y_{t_{n}}=0\quad\mbox{for every }n.\]
\end{lemma}

We also note that this can be extended to a $d$-dimensional setting
where for a $d$-dimensional Wiener process $W(t)$ the zeta process
can be found explicitly \cite{Neate2}.
\begin{theorem}\label{t7}
In $d$-dimensions, the zeta process is given by,
\begin{eqnarray*}
\zeta_t & = & f^0_{(x_t^0(\lambda),t)}(\lambda_1)
-\epsilon x_t^0(\lambda)\cdot W(t)
+\epsilon^2 W(t)\cdot \int_0^t W(s)\di s
+\hf[\epsilon^2]\int_0^t |W(s)|^2 \di s\end{eqnarray*}
where $f^0_{(x,t)}(x_0^1)$ denotes the deterministic reduced action function, $x_t^0(\lambda)$ denotes the pre-parameterisation of the deterministic caustic  and $\lambda$ must satisfy the equation,
\[\nabla_{\lambda}\left(f^0_{(x_t^0(\lambda),t)}(\lambda_1)-
\epsilon x_t^0(\lambda)\cdot W(t)\right)
=0.\]
\end{theorem}
When $\lambda$ is deterministic, the recurrence of this process can be shown using the same argument as for the two dimensional case (further results on recurrence can be found in \cite{Neate4}).
Here we recapitulate our belief that cusped turbulence will be the
most important. As we have shown, when the cusp on the caustic
passes through a level surface, it forces a swallowtail to form on
the level surface. The points of self intersection of this
swallowtail form the Maxwell set.

 \subsection{ Complex turbulence and the resultant $\eta$ process}

We now consider a completely different approach to turbulence.  Let $\left(\lambda,x^2_{0,\mathrm{C}}(\lambda)\right)$ denote the parameterisation of the pre-caustic at time $t$.
When,
\[Z_t=\mbox{Im}\left\{\Phi_t(a+\rmi\eta, x^2_{0,\mathrm{C}}(a+\rmi\eta))\right\},\] is random,
 the values of $\eta(t)$ for which $Z_t=0$ will form a stochastic process.
 The zeros of this new process will correspond to points at which the real pre-caustic touches the complex pre-caustic.
%
The points at which these surfaces touch correspond to swallowtail perestroikas on the caustic.
   When such a perestroika occurs there is a solution of the equations,  \[f'_{({x},t)}(\lambda)=f''_{({x},t)}(\lambda)=
f'''_{({x},t)}(\lambda)=f^{(4)}_{({x},t)}(\lambda)=0.\] Assuming
that $f_{(x,t)}(x_0^1)$ is polynomial in $x_0^1$ we can use the
resultant to state explicit conditions for which this holds
\cite{Neate}.
\begin{lemma}
Let $g$ and $h$ be polynomials of degrees $m$ and $n$ respectively with no common roots or zeros.
Let $f=gh$ be the product polynomial. Then the resultant,
\[ R(f,f')=(-1)^{mn} \left(\frac{m!n!}{N!}\frac{f^{(N)}(0)}{g^{(m)}(0)h^{(n)}(0)}\right)^{N-1}
R(g,g') R(h,h') R(g,h)^2,\]
where $N=m+n$ and $R(g,h)\neq 0$.
\end{lemma}

Since $f'_{(x_t(\lambda),t)}(x_0^1)$ is a polynomial in $x_0$ with real coefficients, its zeros are real or occur in complex conjugate pairs. Of the real roots, $x_0=\lambda$ is repeated. So, \[f'_{(x_t(\lambda),t)}(x_0^1)
= (x_0^1-\lambda)^2 Q_{(\lambda,t)}(x_0^1)
 H_{(\lambda,t)}(x_0^1),\]
 where $Q$ is the product of quadratic factors,
\[Q_{(\lambda,t)}(x_0^1)=\prod\limits_{i=1}^{q}\left\{(x_0^1-a_t^i)^2+(\eta_t^i)^2\right\},\]
 and
 $H_{(\lambda,t)}(x_0^1)$ the product of real factors corresponding to real zeros.
This gives,
\[\left.f'''_{(x_t(\lambda),t)}(x^1_0)\right|_{x^1_0=\lambda}
= 2\prod\limits_{i=1}^q \left\{(\lambda-a_t^i)^2+(\eta_t^i)^2\right\}H_{(\lambda,t)}(\lambda).\]

We now assume that the real roots of $H$ are distinct as are the complex roots of $Q$. Denoting
$\left.f'''_{(x_t(\lambda),t)}(x_0^1)\right|_{x_0^1=\lambda}$ by $f'''_t(\lambda)$ etc, a simple calculation gives
\begin{eqnarray*}
\lefteqn{\left|R_{\lambda}(f'''_t(\lambda),f^{(4)}_t(\lambda))\right|=}\\
 &  &\!\! K_t \prod\limits_{k=1}^q(\eta_t^k)^2\prod\limits_{j\neq k}
\left\{(a_t^k-a_t^j)^4+2((\eta_t^k)^2+(\eta_t^j)^2)(a_t^k-a_t^j)^2+ ((\eta_t^k)^2-(\eta_t^j)^2)^2\right\}\\
& &\quad\times \left|R_{\lambda}(H,H')\right|\left|R_\lambda(Q,H)\right|^2,
\end{eqnarray*}
$K_t$ being a positive constant. Thus, the condition for a swallowtail perestroika to occur is that
\[\rho_{\eta}(t):=\left|R_{\lambda}(f'''_t(\lambda),f^{(4)}_t(\lambda))\right|=0,\]
where we call $\rho_{\eta}(t)$ the \emph{resultant eta process}.

 When the  zeros of $\rho_{\eta}(t)$ form a perfect set, swallowtails will spontaneously
 appear and disappear on the caustic infinitely rapidly.
 As they do so, the geometry of the cool part of the caustic
 will rapidly change as the $\lambda$ shaped sections typical of a swallowtail caustic appear and
 disappear. Moreover, Maxwell sets will be created and destroyed with
 each swallowtail that forms and vanishes adding to the turbulent
 nature of the solution in these regions.
We call this `complex turbulence' occurring at the turbulent times
 which are the zeros of the resultant eta process.

Complex turbulence can be seen as a special case of real turbulence which occurs at specific generalised cusps of the caustic. Recall that when a swallowtail perestroika occurs on a curve, it also satisfies the conditions for having a generalised cusp. Thus, the zeros of the resultant eta process must coincide with some of the zeros of the zeta process for certain forms of cusped turbulence. At points where the complex and real pre-caustic touch, the real pre-caustic and pre-level surface touch in a particular manner (a double touch) since at such a point two swallowtail perestroikas on the level surface have coalesced.

Thus, our separation  of complex turbulence from real turbulence can be seen as an alternative form of categorisation to that outlined in Section 7.1 which could be extended to include other perestroikas.

\end{document}